\documentclass[11pt]{article}
\usepackage{amsmath, amssymb, amsthm, amsbsy, ifthen, graphicx}
\usepackage[active]{srcltx}
\date{}
\title{Packing tight Hamilton cycles in 3-uniform hypergraphs}
\author{
Alan Frieze \thanks{Department of Mathematical
Sciences, Carnegie Mellon University, Pittsburgh, PA 15213, email:
alan@random.math.cmu.edu. Research supported in part by NSF award
DMS-0753472.} 
\and 
Michael Krivelevich \thanks{School of
Mathematical Sciences, Raymond and Beverly Sackler Faculty of Exact
Sciences, Tel Aviv University, Tel Aviv 69978, Israel, e-mail:
krivelev@post.tau.ac.il. Research supported in part by USA-Israel
BSF grant 2006322, by grant 1063/08 from the Israel Science
Foundation, and by a Pazy memorial award.} 
\and
Po-Shen Loh
\thanks{Department of Mathematical Sciences, Carnegie Mellon University, 
Pittsburgh, PA 15217, e-mail: ploh@cmu.edu.} 
}

\usepackage[usenames]{color}\usepackage{graphicx}
\definecolor{brown}{cmyk}{0, 0.72, 1, 0.45}
\definecolor{grey}{gray}{0.5}

\def\e{\epsilon} \def\f{\phi}   
\def\G{\Gamma}

\newtheorem{theorem}{Theorem}[section]
\newtheorem{definition}[theorem]{Definition}
\newtheorem{fact}[theorem]{Fact}
\newtheorem{proposition}[theorem]{Proposition}
\newtheorem{lemma}[theorem]{Lemma}
\newtheorem{corollary}[theorem]{Corollary}

\newcommand{\pr}[1]{\mathbb{P}\left[#1\right]}
\newcommand{\E}[1]{\mathbb{E}\left[#1\right]}

\newcommand{\var}[1]{\text{\rm Var}\left[#1\right]}
\newcommand{\bin}[1]{\text{\rm Bin}\left[#1\right]}
\newcommand{\whp}{\textbf{whp}}
\newcommand{\Gnp}{G_{n,p}}
\newcommand{\andsp}{\quad\quad\quad}
\newcommand{\rt}{\overrightarrow}
\newcommand{\DG}{\rt{\Gamma}}

\newcommand{\DM}{\rt{M}}
\newcommand{\ignore}[1]{}

\def\cE{{\cal E}}
\begin{document}

\maketitle

\begin{abstract}
  Let $H$ be a 3-uniform hypergraph with $n$ vertices.  A tight
  Hamilton cycle $C \subset H$ is a collection of $n$ edges for which
  there is an ordering of the vertices $v_1, \ldots, v_n$ such that
  every triple of consecutive vertices $\{v_i, v_{i+1}, v_{i+2}\}$ is
  an edge of $C$ (indices are considered modulo $n$).  We develop new
  techniques which enable us to prove that under certain natural
  pseudo-random conditions, almost all edges of $H$ can be covered by
  edge-disjoint tight Hamilton cycles, for $n$ divisible by 4.
  Consequently, we derive the corollary that random 3-uniform
  hypergraphs can be almost completely packed with tight Hamilton
  cycles \whp, for $n$ divisible by 4 and $p$ not too small.  Along
  the way, we develop a similar result for packing Hamilton cycles in
  pseudo-random digraphs with even numbers of vertices.
\end{abstract}

\section{Introduction}

Hamilton cycles occupy a position of central importance in graph
theory, and are the subject of countless results.  The most famous is
of course Dirac's Theorem \cite{Dirac}, which states that a Hamilton
cycle can always be found in any $n$-vertex graph with all degrees at
least $n/2$.  Much more work has been done to determine conditions for
Hamiltonicity in graphs, digraphs, hypergraphs, and random and
pseudo-random instances of these objects.  See, e.g., any of
\cite{Bol, HS, KKMR, KrivSud, KMO-type-l, KO-di-survey}.

There has also been a long history of research concerning conditions
for the existence of multiple edge-disjoint Hamilton cycles.  Indeed,
Nash-Williams discovered that the Dirac condition already guarantees
not just one, but at least $\big\lfloor \frac{5}{224} n \big\rfloor$
edge-disjoint Hamilton cycles.  His questions in \cite{NW12, NW13,
  NW14} started a line of investigation, leading to recent work by
Christofides, K\"uhn, and Osthus \cite{CKO}, who answered one of his
conjectures asymptotically by proving that minimum degree $\big(
\frac{1}{2} + o(1) \big) n$ is already enough to guarantee
$\frac{n}{8}$ edge-disjoint Hamilton cycles.  

For random graphs, these ``packings'' with Hamilton cycles are even
more complete.  Bollob\'as and Frieze \cite{BolFrieze} showed that for
every fixed $r$, one can typically find $r$ edge-disjoint Hamilton
cycles in the random graph process as soon as the minimum degree
reaches $2r$.  Kim and Wormald \cite{KW} established a similar result
for random $r$-regular graphs, proving that such graphs typically
contain $\lfloor r/2 \rfloor$ edge-disjoint Hamilton cycles.  The
previous statements are of course best possible, but invite the
natural question of what happens when $r$ is allowed to grow.  Along
these lines, Frieze and Krivelevich showed in \cite{FK-mindegree} that
one can pack $\big\lfloor \frac{\delta}{2} \big\rfloor$ Hamilton
cycles in $\Gnp$, up to $p \leq \frac{(1 + o(1)) \log n}{n}$, where
$\delta$ is the minimum degree of the graph.  For large $p$, they
discovered in \cite{FK0} that one can pack almost all edges into
Hamilton cycles.  This was later improved to essentially the full
range of $p$ by Knox, K\"uhn, and Osthus \cite{KKO}.

In the hypergraph setting, the study of this Hamilton cycle packing
problem was initiated by Frieze and Krivelevich in \cite{FK}.
Although the notion of a Hamilton cycle in an ordinary graph is clear,
there are several ways to generalize the notion to hypergraphs.
Indeed, for any $1 \leq \ell \leq k$, we may define a $k$-uniform
hypergraph $C$ to be a \emph{Hamilton cycle of type $\ell$}\/ if there
is a cyclic ordering of the vertices of $C$ so that every edge
consists of $k$ consecutive vertices, and every pair of consecutive
edges $E_{i-1}, E_i$ in $C$ (according to the natural ordering of the
edges) has $|E_{i-1} \setminus E_i| = \ell$.  The extreme cases $\ell
= 1$ and $\ell = k-1$ are the most obvious generalizations of graph
Hamiltonicity, and cycles of those types are often called
\emph{tight}\/ and \emph{loose}, respectively.  In \cite{FK}, the
first two authors studied the problem of covering almost all the edges
of a given $k$-uniform hypergraph with disjoint Hamilton cycles of a
fixed type $\ell$.  They considered $\ell$ on the looser end of the
spectrum, determining sufficient conditions for the cases $\ell \ge
k/2$.  However, their methods did not extend to the regime $\ell <
k/2$, which seems more difficult.

\subsection{New results}

In this paper, we introduce several new techniques which enable us to
prove the first results for packing \emph{tight}\/ Hamilton cycles,
i.e., with $\ell = 1$.  To give the clearest presentation of the new
methods involved, we concentrate on the concrete case of 3-uniform
hypergraphs, which we refer to as \emph{3-graphs}\/ from now on.  Our
approach can be naturally extended to the general $k$-uniform case,
but the explanations necessarily become more involved.

We focus on 3-graphs with certain pseudo-random properties.
Consequently, our results will yield corollaries in the standard
random 3-graph model, denoted $H_{n,p;3}$, in which each of the
$\binom{n}{3}$ possible edges independently appears with probability
$p$.  Recall that there are several essentially equivalent notions of
pseudo-randomness in graphs, most notably the ones involving
uniformity of edge distribution, the second eigenvalue, or the global
count of 4-cycles.  (See, e.g., the survey \cite{KS-survey}.)
However, even in ordinary graphs, these global conditions are
insufficient for Hamilton cycle packing, because Hamilton cycles visit
every vertex.  Therefore, even a single non-conforming vertex can doom
the entire packing process.

This situation can be prevented by localizing the specification of the
pseudo-random criteria.  Importantly, the resulting stronger
conditions are still satisfied by the random objects in question.  For
example, in the graph case, this can be done by controlling the
degrees and codegrees (see, e.g., \cite{CGW, Thomason}).  Note that
this is a special case of the more general notion of controlling the
numbers of \emph{extensions}\/ to certain specific subgraphs,
uniformly over all base sites, instead of just the global number of
copies of a certain graph ($C_4$, say).  The following definition
generalizes that concept to 3-graphs.

\begin{definition}
  \label{def:3graph-ext} 
  Let $v_1, \ldots, v_t$ be distinct vertices of a 3-graph $H$, and
  let $\Gamma$ be a graph on vertex set $[t]$.  Then, we define
  $d_\Gamma(v_1, \ldots, v_t)$ to be the number of vertices $x \in H$
  such that $\{v_i, v_j, x\} \in H$ for every edge $ij \in \Gamma$.
\end{definition}

\noindent \textbf{Example 1.}\, If $t = 2$ and $\Gamma$ is a single
edge, then $d_\Gamma(v_1, v_2)$ is the number of edges which contain
the pair $\{v_1, v_2\}$.  This can be interpreted as a generalized
notion of degree.

\vspace{3mm}

\noindent \textbf{Example 2.}\, In general, if $\Gamma$ is the
complete graph on $t$ vertices, then $d_\Gamma (v_1, \ldots, v_t)$ is
the number of $x$ which simultaneously form edges with all pairs
$\{v_i, v_j\}$.  This can be interpreted as a generalized notion of
codegree.

\vspace{3mm}

We can now state our localized version of pseudo-randomness, which we
will later show is satisfied by the random 3-graph $H_{n,p;3}$
\whp\footnote{Here a sequence of events $\cE_n,n\geq 0$ is said to
  occur {\em with high probability} (abbreviated \whp), if
  $\lim_{n\to\infty}\Pr(\cE_n)=1$.} when $p$ is not too small.  Since
we need to control many quantities within certain ranges, we will
employ the notation $A=(1 \pm \epsilon)B$ as a shorthand for the pair
of inequalities $(1 - \epsilon)B \leq A \leq (1 + \epsilon)B$.

\begin{definition}
  \label{def:3graph-props} We say that an $n$-vertex 3-graph $H$ is
  \textbf{$\boldsymbol{(\epsilon, p)}$-uniform} if for every auxiliary
  graph $\Gamma$ on $t \leq 7$ vertices and $s \leq 6$ edges, and
  every choice of distinct vertices $v_1, \ldots, v_t \in H$, we have
  \begin{displaymath}
    d_\Gamma (v_1, \ldots, v_t) = (1 \pm \epsilon) np^s.
  \end{displaymath}
\end{definition}

\noindent \textbf{Remark.}\, We do not actually need the full strength
of this definition, as will be clear in the proof of Lemma \ref{lem:d}
in Section \ref{sec:3graph}.  Indeed, there are only 5 specific graphs
$\Gamma$ with respect to which we require control (depicted in Figure
\ref{fig:gamma} of Section \ref{sec:3graph}).  However, we feel that
the description is more succinctly captured in the above statement.

\vspace{3mm}

The main result of this paper establishes the first known packing
result for tight Hamilton cycles in pseudo-random hypergraphs.  Here,
and in the rest of this paper, we write $f(n) \ll g(n)$ if $f/g
\rightarrow 0$ as $n \rightarrow \infty$.

\begin{theorem}
  \label{thm:3graph}
  Suppose that $n$ is a sufficiently large multiple of four, and
  $\epsilon^{45} n p^{16} \gg \log^{21} n$.  Then every $(\epsilon,
  p)$-uniform 3-graph with $n$ vertices can have all but at most
  $\epsilon^{1/15}$-fraction of its edges covered by a disjoint union
  of tight Hamilton cycles.
\end{theorem}

Note that in $H_{n,p;3}$, for each graph $\Gamma$ with $s \leq 6$
edges on a set of $t\leq 7$ vertices, the value of $d_\G$ has
distribution $\bin{n-t,p^s}$. So, the Chernoff bound in Section
\ref{sec:tools} (Fact \ref{fact:chernoff}) shows that as long as
$\epsilon^2 n p^6 \gg \log n$, we have
\begin{displaymath}
  \pr{H_{n,p;3} \text{ is not $(\e,p)$-uniform} }
  =
  O(n^7) \cdot
  \sum_{s=1}^6 \pr{\bin{n,p^s} \neq (1 \pm \e) np^s}
  =
  o(1),
\end{displaymath}
giving the following immediate corollary.

\begin{corollary}
  \label{cor:3graph}
  Suppose that $\epsilon, n, p$ satisfy $\epsilon^{45} n p^{16} \gg
  \log^{21} n$.  Then whenever $n$ is a multiple of four, $H_{n,p;3}$
  can have all but at most $\epsilon^{1/15}$-fraction of its edges
  covered by a disjoint union of tight Hamilton cycles \whp.
\end{corollary}

\noindent \textbf{Remark.}\, Although both results are stated for $n$
divisible by 4, we expect that they are true in general.  Note,
however, that a divisibility condition is unavoidable in the general
case of packing Hamilton cycles of type $\ell$ in $k$-uniform
hypergraphs, since $\ell$ must divide $n$.

\vspace{3mm}

Along the way, we also prove a new result about packing Hamilton
cycles in pseudo-random digraphs.  The result differs from that in
\cite{FK0} as our definition of pseudo-randomness is local instead of
global, and therefore may be easier to apply in some situations.
Indeed, the previous result required a minimum degree condition,
together with bounds on the edge distributions across almost all cuts.
We can now replace the latter conditions with two more local
statements.  The specific conditions that we impose below have been
tailored for the task of producing Hamilton cycles in digraphs.

\begin{definition}
  \label{def:digraph-props} We say that an $n$-vertex digraph is
  \textbf{$\boldsymbol{(\epsilon, p)}$-uniform} if it satisfies the
  following properties:
  \begin{description}
  \item[(i)] Every vertex $a$ has out-degree $d^+(a) = (1 \pm \epsilon)
    np$ and in-degree $d^-(a) = (1 \pm \epsilon) np$.

  \item[(ii)] For every pair of distinct vertices $a, b$, all three of
    the following quantities are $(1 \pm \epsilon) np^2$: the number of
    common out-neighbors $d^+(a, b)$, the number of common
    in-neighbors $d^-(a, b)$, and the number of out-neighbors of $a$
    which are also in-neighbors of $b$.

  \item[(iii)] Given any four vertices $a, b, c, d$, which are all
    distinct except for the possibility $b=c$, there are $(1 \pm
    \epsilon) np^4$ vertices $x$ such that $\rt{ax}$,
    $\rt{xb}$, $\rt{cx}$,
    $\rt{xd}$ are all directed edges.
  \end{description}
\end{definition}

\noindent \textbf{Remark.}\, It is not clear that this is the minimal
set of pseudo-random conditions which enable Hamilton cycle packing in
directed graphs.  We choose the above statements because they
naturally arise from our analysis, and are therefore the most
convenient for our purposes.  Importantly, one can easily see that
they are satisfied by random digraphs \whp, as long as the edge
probability is not too small.

\vspace{3mm}

Under these easily-verifiable conditions, we are able to prove the
following packing result for digraphs, which has the obvious corollary
for random digraphs which are not too sparse.

\begin{theorem}
  \label{thm:digraph}
  Suppose that $\epsilon^{11} n p^8 \gg \log^5 n$, and $n$ is a
  sufficiently large even integer.  Then every $(\epsilon, p)$-uniform
  digraph can have its edges partitioned into a disjoint union of
  directed Hamilton cycles, except for a set of at most
  $\epsilon^{1/8}$-fraction of its edges.
\end{theorem}

Although we originally developed this result only as a building block
for our 3-graph analysis in Section \ref{sec:3graph}, we feel it is
worth bringing attention to, as it may be of independent interest.  In
particular, it is easier to apply than its counterpart in \cite{FK0},
because our pseudo-randomness conditions are easier to verify.

\subsection{Proof overview and organization}
\label{sec:overview}

The key insight in the proof of Theorem \ref{thm:3graph} is the
following connection between tight Hamilton cycles in $H$ and Hamilton
cycles in an associated digraph.  For a random permutation $v_1, v_2,
\ldots, v_n$ of the vertices of $H$, define an $\frac{n}{2}$-vertex
digraph $D$ with vertex set $\{(v_1, v_2), (v_3, v_4), \ldots,
(v_{n-1}, v_n)\}$.  Note that each vertex of $D$ corresponds to an
ordered pair of vertices of $H$, so $D$ will have an even number of
vertices, since the number of vertices of $H$ is a multiple of 4.
Place a directed edge from $(v_i, v_{i+1})$ to $(v_j, v_{j+1})$ if and
only if both hyperedges $\{v_i, v_{i+1}, v_j\}$ and $\{v_{i+1}, v_j,
v_{j+1}\}$ are present in $H$.  In this construction, Hamilton cycles
in $D$ give rise to tight Hamilton cycles in $H$.

To extract edge disjoint Hamilton cycles from a digraph $D$ with an
even number of vertices, we use an approach similar to that taken in
\cite{FK}.  Let $w_1, w_2, \ldots, w_{2m}$ be a random permutation of
the vertices of $D$ with $m = n/4$, and define $A = \{w_1, w_2,
\ldots, w_m\}$ and $B=\{w_{m+1}, \ldots, w_{2m}\}$.  Define a
bipartite graph $\G$ with bipartition $(A,B)$, and place an edge
between $w_i \in A$ and $w_j \in B$ whenever $\rt{w_i w_j}$ and
$\rt{w_j w_{i+1}}$ are both edges of $D$.  Now perfect matchings in
$\Gamma$ give rise to Hamilton cycles in $D$, and previous approaches
in \cite{FK} show how to pack perfect matchings in pseudo-random
bipartite graphs.

However, not all Hamilton cycles in $D$ arise from perfect matchings
in one particular $\Gamma$.  Similarly, not all Hamilton cycles in $H$
arise from Hamilton cycles in a single $D$.  We overcome both
obstacles with the same iterative approach, which we illustrate for
the hypergraph packing.  Roughly speaking, instead of stopping after
generating a single $D$, we sequentially generate digraphs $D_1, D_2,
\ldots, D_r$ in the above manner, extracting a large set of edge
disjoint directed Hamilton cycles from each, and deleting the
corresponding edge-disjoint Hamilton cycles from $H$.  At each step,
we verify that the pseudo-random properties are maintained.  We repeat
the process until we have packed the required number of cycles.

Since the digraph packing and hypergraph packing proofs are
essentially independent (though similar), we separate them into
Sections \ref{sec:digraph} and \ref{sec:3graph}, respectively.  All
required concentration inequalities are collected in Section
\ref{sec:tools} for the reader's convenience.  The final section
contains some concluding remarks and open problems.

\subsection{Notation and conventions}

We will implicitly assume throughout that $\epsilon, p$ are small,
e.g., less than $1/10$.  Our results have them tending to zero.  The
following (standard) asymptotic notation will be utilized extensively.
For two functions $f(n)$ and $g(n)$, we write $f(n) = o(g(n))$ or
$g(n) = \omega(f(n))$ if $\lim_{n \rightarrow \infty} f(n)/g(n) = 0$,
and $f(n) = O(g(n))$ or $g(n) = \Omega(f(n))$ if there exists a
constant $M$ such that $|f(n)| \leq M|g(n)|$ for all sufficiently
large $n$.  We also write $f(n) = \Theta(g(n))$ if both $f(n) =
O(g(n))$ and $f(n) = \Omega(g(n))$ are satisfied.
All logarithms will be in base $e \approx 2.718$.

\section{Probabilistic tools}
\label{sec:tools}

We recall the Chernoff bound for exponential concentration of the
binomial distribution.  The following formulation appears in, e.g.,
\cite{AS}.

\begin{fact}
  \label{fact:chernoff}
  For any $\epsilon > 0$, there exists $c_\epsilon > 0$ such that any
  binomial random variable $X$ with mean $\mu$ satisfies
  \begin{displaymath}
    \pr{|X - \mu| > \epsilon \mu} < e^{-c_\epsilon \mu},
  \end{displaymath}
  where $c_\epsilon$ is a constant determined by $\epsilon$.  When
  $\epsilon < 1$, we may take $c_\epsilon = \frac{\epsilon^2}{3}$.
\end{fact}

A binomial random variable is the sum of independent indicator
variables.  We also need concentration in settings which are still
product spaces, but are somewhat more complicated than simple
binomials.  A random variable $X(\omega)$ defined over an
$n$-dimensional product space $\Omega = \prod_{i=1}^n \Omega_i$ is
called $C$-Lipschitz if changing $\omega$ in any single coordinate
affects the value of $X(\omega)$ by at most $C$.  The Hoeffding-Azuma
inequality (see, e.g., \cite{AS}) provides concentration for
these distributions.

\begin{fact}
  \label{fact:azuma}
  Let $X$ be a $C$-Lipschitz random variable on an $n$-dimensional
  product space.  Then for any $t \geq 0$,
  \begin{displaymath}
    \pr{ | X - \E{X} | > t }
    \leq
    2 \exp\left\{
      -\frac{t^2}{2 C^2 n}
    \right\}.
  \end{displaymath}
\end{fact}

We also need concentration in settings where the probability space is
not a simple $n$-dimensional product space, but rather the set of
permutations on $n$ elements.  The following concentration inequality
can be found in, e.g., Section 3.2 of \cite{McDiarmid} or Lemma 11 of
\cite{FP}.

\begin{fact}
  \label{fact:conc-perm}
  Let $X$ be a random variable on the uniformly distributed
  probability space of permutations on $n$ elements, and let $C$ be a
  real number.  Suppose that whenever $\sigma, \sigma' \in S_n$ differ
  by a single transposition, $|X(\sigma) - X(\sigma')| \leq C$.  Then,
  \begin{displaymath}
    \pr{\left| X - \E{X} \right| \geq t} 
    \leq 
    2\exp\left\{  
      -\frac{2t^2}{C^2 n}
    \right\}.
  \end{displaymath}
\end{fact}

\section{Packing Hamilton cycles in digraphs}
\label{sec:digraph}

In this section, we show how to complete the outline in Section
\ref{sec:overview} for Hamilton cycle packing in digraphs.  Recall
that the main idea in this part is to reduce the problem to packing
perfect matchings in bipartite graphs.  We begin by showing how to
achieve the final step.

\subsection{Packing perfect matchings}
\label{sec:graph}

Here, we will show that we can efficiently decompose a suitably
uniform bipartite graph into perfect matchings.  This is contained in
the following proposition, whose statement is very similar to Lemma 6
of \cite{FK}.  For completeness, we provide a slightly different proof
here.


\begin{proposition}
  \label{prop:graph}
  Let $G$ be a bipartite graph with parts $A$ and $B$, each of size
  $m$.  Suppose that $p, \epsilon < 1/2$ are given such that every
  vertex has degree $(1 \pm \epsilon)mp$, and every pair of distinct
  vertices has codegree at most $(1+\epsilon) mp^2$.  Also assume that
  $\epsilon^{4/3} mp^2 \geq 1$.  Then the edges of $G$ can be
  partitioned into the disjoint union $E_0 \cup \ldots \cup E_k$,
  where $|E_0| < 4 \epsilon^{1/3} e(G)$ and each $E_i$ with $i \geq 1$
  is a perfect matching.
\end{proposition}

We mentioned in the introduction that although pseudo-randomness can
be deduced from the global 4-cycle count, local conditions are
required to pack spanning objects such as Hamilton cycles or perfect
matchings.  For bipartite graphs, the above proposition's conditions
on codegrees and degrees provide this localized property.  Indeed, the
concentrated codegree condition implies the global bound on the number
of 4-cycles, simply by summing over all vertex pairs.  The following
lemma converts this into a suitable bound on the uniformity of edge
distribution---another pseudo-random property, which we actually need
to complete this proof.


\begin{lemma}
  \label{lem:flow-prelim}
  Let $G$ be a bipartite graph with parts $A$ and $B$, each of size
  $m$.  Let $X$ and $Y$ be subsets of $A$ and $B$, respectively, with
  $|X| \geq 1/(\epsilon p)$ and $|Y| \geq \epsilon^{1/3} m$.  Suppose
  that $\epsilon, p$ are given such that every vertex in $X$ has
  degree $(1 \pm \epsilon)mp$, and every pair of distinct vertices in
  $X$ has codegree at most $(1+\epsilon) mp^2$.  Then $e(X, Y) \geq (1
  - 3\epsilon^{1/3}) |X||Y|p$.
\end{lemma}

\noindent \textbf{Proof.}\, Let $x = |X|$ and $y = |Y|$.  For each
vertex $v \in B$, let $d_X(v)$ denote the number of neighbors $v$ has
in $X$.  Define the random variable $Z$ by sampling a uniformly random
vertex $v \in B$, and setting $Z = d_X(v)$.  Since every vertex in $X$
has degree at least $(1-\epsilon)mp$, we have $e(X, B) \geq
(1-\epsilon) xmp$, and so $\E{Z} \geq (1-\epsilon) xp$.

{From} the codegree condition, the number of labeled paths of length 2
from $X$ to $X$ is at most $(1+\epsilon) x^2 mp^2$.  This quantity is
also equal to 
\begin{displaymath}
  \sum_{v \in B} d_X(v) \cdot (d_X(v) - 1)
  =
  \sum_{v \in B} d_X(v)^2  - e(X, B),
\end{displaymath}
so $\sum_{v \in B} d_X(v)^2 \leq (1+\epsilon) x^2 mp^2 + (1+\epsilon)
xmp$.  Therefore, $\E{Z^2} \leq (1+\epsilon) (x^2 p^2 + xp)$, which is
at most $(1+3\epsilon)(xp)^2$ since we assumed $x \geq 1/(\epsilon p)$.

This implies that $\var{Z} = \E{Z^2} - \E{Z}^2 \leq 5\epsilon(xp)^2$,
and this low variance will allow us to conclude that vertices in $Y$
cannot have an average $d_X$ which is too low.  Formally, if we assume
for contradiction that $e(X, Y) < (1 - 3\epsilon^{1/3}) xyp$, then
Jensen's inequality gives (here, $\overline{Z}$ is a shorthand for
$\E{Z}$):
\begin{eqnarray*}
  \E{(Z-\overline{Z})^2 \mid v \in Y}
  &\geq& 
  \left( \E{Z-\overline{Z} \mid v \in Y} \right)^2 \\
  &=&
  \left( \frac{e(X, Y)}{|Y|} - \frac{e(X, B)}{|B|} \right)^2 \\
  &\geq&
  \left( (3\epsilon^{1/3} - \epsilon)(xp) \right)^2 \\
  &>& 5 \epsilon^{2/3} (xp)^2.
\end{eqnarray*}
Yet a uniformly random vertex in $B$ lies in $Y$ with probability at
least $\epsilon^{1/3}$, so
\begin{displaymath}
  \var{Z} 
  \geq 
  \E{(Z-\overline{Z})^2 \mid v \in Y} \cdot \pr{v \in Y} 
  > 
  5\epsilon (xp)^2,
\end{displaymath}
contradiction.  \hfill $\Box$

\vspace{3mm}

Now that we know the edges are distributed fairly uniformly, we can
prove the packing result using the maxflow-mincut theorem.

\vspace{3mm}

\noindent \textbf{Proof of Proposition \ref{prop:graph}.}\, First,
observe that if we can pack $k = (1 - 3\epsilon^{1/3}) mp$
edge-disjoint perfect matchings into $G$, then the proportion of
unused edges is at most
\begin{displaymath}
  \frac{\epsilon + 3\epsilon^{1/3}}{1+\epsilon} \leq 4 \epsilon^{1/3}.
\end{displaymath}
We will use the maxflow-mincut theorem to show that we can indeed pack
$k$ matchings.  Observe that one can pack $k$ edge-disjoint perfect
matchings in $G$ if and only if the following network has an integer
flow of size exactly $km$.  Give each edge in $G$ a capacity of 1, and
direct it from $A$ to $B$.  Add a source vertex $\sigma$, with an edge
of capacity $k$ to each vertex in $A$, and add a sink vertex $\tau$,
with an edge of capacity $k$ from each vertex in $B$.  Clearly, an
integer flow of size $km$ is achieved if and only if the subset of
used edges between $A$ and $B$ forms a $k$-regular graph.  Since every
$k$-regular bipartite graph can be decomposed into $k$ edge-disjoint
perfect matchings, this completes the argument.

All capacities are integers, so there is an integer flow which
achieves the maximum size.  Therefore, by the maxflow-mincut theorem,
it remains to show that every cut has size at least $km$.  Consider an
arbitrary cut.  Let $X \subset A$ be those vertices of $A$ which are
on the same side as $\sigma$, and let $Y \subset B$ be those vertices
of $B$ which are on the same side as $\tau$.  Let $x = |X|$ and $y =
|Y|$.  The size of this cut is then exactly $k(m-x) + k(m-y) + e(X,
Y)$.  Therefore, it suffices to establish the following inequality for
all choices of $X \subset A$ and $Y \subset B$:
\begin{equation}
  \label{ineq:cut}
  k(m-x) + k(m-y) + e(X, Y) \geq km
\end{equation}
This is purely an inequality about the original graph $G$.  Since it
is symmetric in $X$ and $Y$, assume without loss of generality that $x
\leq y$.

We will use Lemma \ref{lem:flow-prelim} to estimate $e(X, Y)$, but
first, we will need to dispose of the boundary cases $y \leq
\epsilon^{1/3} m$ and $x \leq 1/(\epsilon p)$, in which it does not
apply.  In the first case, observe that we automatically have $x \leq
y \leq \epsilon^{1/3} m$ as well, so $k(m-x) + k(m-y)$ is already at
least $km$.  (We may assume $\epsilon^{1/3} < 1/2$.)  In the second
case, note that if $m-y \geq x$, then $k(m-x) + k(m-y)$ is already at
least $km$.  So, we may assume that $m-y \leq x \leq 1/(\epsilon p)$.
Yet every vertex in $X$ has at least $(1-\epsilon)mp$ neighbors in
$B$, and only at most $m-y$ of them can be outside $Y$.  Therefore,
\begin{displaymath}
  e(X, Y) 
  \geq 
  x \cdot [ (1-\epsilon)mp - (m-y) ]
  \geq 
  x \cdot \left[ (1-\epsilon)mp - \frac{1}{\epsilon p} \right]
  \geq
  x \cdot (1 - 2\epsilon^{1/3})mp,
\end{displaymath}
since we assumed that $1/(\epsilon p) \leq \epsilon^{1/3} mp$.  Yet
the final quantity exceeds $xk$, so we also have inequality
\eqref{ineq:cut} in this case.

In all remaining cases, we see that Lemma \ref{lem:flow-prelim} applies, so
\begin{displaymath}
  k(m-x) + k(m-y) + e(X, Y)
  \geq
  k(m-x) + k(m-y) + (1-3\epsilon^{1/3})xyp.
\end{displaymath}
The right-hand side is bilinear in $x$ and $y$, so it is minimized
when $x, y \in \{0, m\}$.  If either of $x$ or $y$ are 0, then one of
the first two terms is already $km$.  On the other hand, if $x = y =
m$, then the expression is precisely $(1-3\epsilon^{1/3}) m^2 p = km$ as
well.

Therefore, all cuts in the network have size at least $km$, so by the
maxflow-mincut theorem, there is a flow of size $km$, completing the
proof. \hfill $\Box$

\subsection{Directed Hamilton cycles and perfect matchings}

Now that Proposition \ref{prop:graph} allows us to efficiently pack
perfect matchings, we can use it as the base to which we reduce the
problem of packing Hamilton cycles in directed graphs.  Throughout
this section, all directed graphs have no loops, and no repeated edges
in the same direction.  However, they may have edges in both
directions between pairs of vertices.

The main objective of this section is to prove Theorem
\ref{thm:digraph}, which allows us to efficiently pack Hamilton cycles
in digraphs that have the pseudo-random properties specified in
Definition \ref{def:digraph-props}.  As mentioned in the proof
overview in Section \ref{sec:overview}, we extract Hamilton cycles
from digraphs by connecting them to perfect matchings in ordinary
graphs.  Consider the following random procedure, which constructs an
undirected bipartite graph $\Gamma$ from a digraph $D$ with an even
number of vertices.

\vspace{3mm}

\noindent \textbf{Procedure 1.}\, This takes as input a digraph $D$
with an even number of vertices.
\begin{enumerate}
\item Generate a random permutation $\sigma = (v_1, v_2, \ldots, v_n)$
  of the vertices of $D$.  Consider this sequence as two consecutive
  segments of length $\frac{n}{2}$, and let $A = \{v_1, \ldots,
  v_{\frac{n}{2}}\}$ and $B = \{v_{\frac{n}{2} + 1}, \ldots, v_n\}$.
  The graph $\Gamma$ will be a bipartite graph with parts $A$ and $B$.

\item Define \emph{the successor of $v_i$}\/ to be $v_{i+1}$, unless
  $i = \frac{n}{2}$ (in which case the successor is $v_1$) or $i = n$
  (in which case it is $v_{\frac{n}{2} + 1}$).  Similarly, define
  \emph{the predecessor of $v_i$}\/ to be $v_{i-1}$, unless $i=1$ (in
  which case it is $v_{\frac{n}{2}}$), or $i = \frac{n}{2} + 1$ (in
  which case it is $v_n$).

\item For $v_i \in A$ and $v_j \in B$, place the edge $v_i v_j$ in
  $\Gamma$ if and only if $D$ contains both directed edges $\rt{v_i
    v_j}$ and $\rt{v_j v_{i^+}}$, where $v_{i^+}$ is the successor of
  $v_i$ as defined in the previous step.

\item Since an edge in $\Gamma$ corresponds to two edges in $D$, we
  account for this by defining $\DG \subset D$ to be the digraph
  containing all $\{\rt{v_i v_j}, \rt{v_j v_{i^+}}\}$, for each $v_i
  v_j \in \Gamma$ with $v_i \in A$.
\end{enumerate}

\vspace{3mm}

The key observation is that perfect matchings in $\Gamma$ cleanly
correspond to Hamilton cycles in $\DG$.  Specifically, given a perfect
matching $M$ in $\Gamma$, one can recover a Hamilton cycle in $\DG$.
Indeed, for each $v_i \in A$, $M$ matches $v_i$ to a distinct $v_{i'}
\in B$.  So, $(v_1, v_{1'}, v_2, v_{2'}, \ldots, v_{n/2}, v_{(n/2)'})$
is a directed Hamilton cycle in $\DG$; call this the Hamilton cycle
\emph{associated with}\/ $M$.

\begin{lemma}
  \label{lem:hamilton-matching}
  Let $M_1, M_2$ be edge-disjoint perfect matchings in $\Gamma$.  Then
  their associated Hamilton cycles $\DM_1, \DM_2$ are also
  edge-disjoint in $\DG$.
\end{lemma}

\noindent \textbf{Proof.}\, Since $M_1$ and $M_2$ are disjoint, it is
clear that $\DM_1$ and $\DM_2$ cannot overlap on any edges directed
from $A$ to $B$.  On the other hand, if both $\DM_i$ contain the same
edge $\rt{v_j v_k}$, where $v_j \in B$ and $v_k \in A$,
then they also both contain $\rt{v_{k-1} v_j}$, by
definition of $\DM_i$.  Then both $M_i$ contain the edge $v_{k-1}
v_j$, contradiction.  \hfill $\Box$

\vspace{3mm}

The next step, as we mentioned in the proof overview in the
introduction, is to show that the randomly constructed bipartite graph
inherits the pseudo-random properties of the initial digraph.

\begin{lemma}
  \label{lem:gamma}
  Suppose $n$, $p$, and $\epsilon$ satisfy $\epsilon^2 np^8 \gg \log
  n$.  Let $D$ be an $(\epsilon, p)$-uniform digraph on $n$ vertices
  ($n$ even), and randomly construct the undirected bipartite graph
  $\Gamma$ according to Procedure 1.  Then, with probability
  $1-o(n^{-1})$, $\Gamma$ satisfies all of the following properties:
  \begin{description}
  \item[(i)] All degrees are $(1 \pm 5\epsilon) \frac{np^2}{2}$.
  \item[(ii)] All codegrees between pairs of vertices on the same side
    of the bipartition are $(1 \pm 5\epsilon) \frac{np^4}{2}$.
  \end{description}
\end{lemma}

\noindent \textbf{Proof.}\, By construction, $\Gamma$ is a bipartite
graph with parts $A = \{v_1, \ldots, v_{\frac{n}{2}}\}$ and $B =
\{v_{\frac{n}{2} + 1}, \ldots, v_n\}$.  There are essentially four
claims to prove, as the parts $A$ and $B$ are not constructed
symmetrically.

\vspace{3mm}

\noindent \textbf{Claim 1: Degrees in $\boldsymbol{A}$ are correct.}\,
To show this, for an arbitrary vertex $a$, define a random variable
$N_a$ as follows.  Let $b$ be the successor of $a$, 
as defined in Step 2 of Procedure 1.  $N_a$ is of no interest when $a\in B$ but allowing $a\in B$ makes the proof a trifle simpler.
By definition, $a$ and $b$ are either both in $A$ or
both in $B$; let $N_a$ be the number of vertices $x$ in the other part
such that both $\rt{ax}$ and $\rt{xb}$ are
edges of $D$.  Note that $N_a$ is actually the degree of $a$ in
$\Gamma$ if $a$ happens to end up in $A$.  Therefore, it suffices to
show that $N_a = (1 \pm 3\epsilon) \frac{np^2}{2}$ with probability $1
- o(n^{-2})$.

Note that the random variable $N_a$ is completely determined by the
permutation $\sigma$.  Expose $\sigma$ in stages.  First, expose which
index $i$ has $v_i = a$.  Next, expose the identity of the vertex $b$
which is the successor of $a$.  This determines $v_i$ and $v_{i^+}$,
where $i^+$ is either $i+1$, $1$, or $n/2 + 1$.  Conditioned on these,
the associations to the remaining vertices $v_j$ are a uniform
permutation over the remaining $n-2$ vertices.  Yet the last part of
property (ii) for $(\epsilon, p)$-uniformity shows that in $D$, the
number of vertices $x$ with $\rt{ax}$ and
$\rt{xb}$ in $D$ is $(1 \pm \epsilon)np^2$.  As each $x$ will
land in the part opposite $\{a, b\}$ with probability exactly
$\frac{n/2}{n - 2}$, we immediately have $\E{N_a} = (1 \pm
2\epsilon)\frac{np^2}{2}$.

For concentration, we use Fact \ref{fact:conc-perm}.  We already
conditioned on $a$ and $b$, so the remaining randomness is from a
uniformly random permutation of $n-2$ elements.  Transposing two of
those elements can only change $N_a$ by at most 1, so by Fact
\ref{fact:conc-perm}, $N_a$ deviates from its expectation by $\epsilon
\cdot \frac{np^2}{2}$ with probability at most
\begin{displaymath}
  2\exp\left\{-2 \left(\epsilon
      \cdot \frac{np^2}{2} \right)^2 \bigg/ (n-2)
  \right\}
  \leq
  \exp\left\{- \frac{ \epsilon^2 n p^4 }{3} \right\}
  = 
  o(n^{-2}).
\end{displaymath}
Therefore, $N_a$ is indeed $(1 \pm 3\epsilon)\frac{np^2}{2}$ with
probability $1 - o(n^{-2})$, as desired.

\vspace{3mm}

\noindent \textbf{Claim 2: Degrees in $\boldsymbol{B}$ are correct.}\,
Consider an arbitrary vertex $b$.  Let $S_b$ be the set of ordered
pairs of distinct vertices $(x, y)$ with $\rt{xb},
\rt{by} \in D$.  By property (i) of $(\epsilon,
p)$-uniformity, $b$ has $(1 \pm \epsilon) np$ in-neighbors, and $(1
\pm \epsilon) np$ out-neighbors.  Since $\epsilon np \gg 1$, 
this
implies that $|S_b| = (1 \pm 3\epsilon) n^2 p^2$.  Let $N_b$ be the
number of these pairs that, after the permutation $\sigma$, have the
additional properties that \textbf{(i)} both $x$ and $y$ are on the
opposite side of the bipartition to $b$, and \textbf{(ii)} $y$ is
the successor of $x$.  It suffices to show that $N_b$ is $(1 \pm
5\epsilon) \frac{np^2}{2}$ with probability $1 - o(n^{-2})$, because
$N_b$ is precisely the degree of $b$ whenever $b$ happens to end up in
$B$.

For this, we employ the same strategy as used in the previous claim.
First expose the index $i$ for which $v_i = b$.  Then, the locations
of the remaining vertices are uniformly permuted amongst the remaining
$n-1$ positions.  For a given ordered pair $(x, y) \in S_b$, the
probability that $x$ lands on the opposite side of $b$ is precisely
$\frac{n/2}{n-1}$.  Then, conditioned on this, the probability that
$y$ is the successor of $x$ is precisely $\frac{1}{n-2}$.  Therefore,
\begin{displaymath}
  \E{N_b} 
  = 
  |S_b| \cdot \frac{n/2}{n-1} \cdot \frac{1}{n-2}
  =
  (1 \pm 4\epsilon) \frac{np^2}{2}.
\end{displaymath}
For concentration, we again expose $b$ first, and then consider the
resulting $(n-1)$-permutation.  Consider two such permutations
$\sigma$ and $\sigma'$ differing only on a single transposition.
Transposing two of those elements can only change $N_b$ by at most 4,
so by Fact \ref{fact:conc-perm}, the probability that $N_b$ deviates
by over $\epsilon \cdot \frac{np^2}{2}$ from its expectation is at
most
\begin{displaymath} 
  2\exp\left\{-
    2 \left(
      \epsilon \cdot \frac{np^2}{2}
    \right)^2 / ({4}^2 (n-1))
  \right\}
  =
  o(n^{-2}),
\end{displaymath}
as before.  Hence $N_b$ is indeed $(1 \pm 5\epsilon) \frac{np^2}{2}$
with the desired probability $1 - o(n^{-2})$.

\vspace{3mm}

\noindent \textbf{Claim 3: Co-degrees in $\boldsymbol{A}$ are
  correct.}\, This is similar to Claim 1.  Fix any two distinct
vertices $a$ and $b$, and let $c$ and $d$ be their respective
successors as defined in Step 1 of
Procedure 1.  Note that we may have $b=c$ or $a=d$, but not both.  Let
$N_{a,b}$ be the number of vertices $x$ such that $x$ is on the
opposite side of $a$, and $\rt{ax}, \rt{xc},
\rt{bx}, \rt{xd} \in D$.

Expose the indices $i, j$ for which $v_i = a$ and $v_j = b$, and
expose their successors $c$ and $d$.  The remaining vertices are
uniformly permuted over the remaining $n-4$ positions (or $n-3$ positions if $b = c$ or $a = d$, although the argument will be the same).  By property
(iii) of $(\epsilon, p)$-uniformity, there are $(1 \pm \epsilon) np^4$
candidates for $x$.
{As each $x$ will
land in the part opposite $a$ with probability either
$\frac{n/2}{n - 4}$ (if $b$ is on the same side as $a$)
or $\frac{n/2 - 2}{n - 4}$ (if $b$ is opposite $a$),}
it follows
that $\E{N_{a,b}} = (1 \pm 2\epsilon) \frac{np^4}{2}$.  A
transposition in the $(n-4)$-permutation can only affect $N_{a,b}$ by
at most 1, so the probability that $N_{a,b}$ deviates from its
expectation by over $\epsilon \cdot \frac{np^4}{2}$ is at most
\begin{displaymath} 
  2\exp\left\{-
    2 \left(
      \epsilon \cdot \frac{np^4}{2}
    \right)^2 / (n-4)
  \right\}
  =
  o(n^{-3}).
\end{displaymath}
Taking a union bound over all pairs $(a, b)$ yields the desired
result.

\vspace{3mm}

\noindent \textbf{Claim 4: Co-degrees in $\boldsymbol{B}$ are
  correct.}\, This is similar to Claim 2.  Fix any two distinct
vertices $a$ and $b$, and let $S_{a,b}$ be the set of ordered pairs of
distinct vertices $(x, y)$ with $\rt{xa}, \rt{ay}, \rt{xb}, \rt{by}
\in D$.  By property (ii) of $(\epsilon, p)$-uniformity, $a$ and $b$
have $(1 \pm \epsilon) np^2$ common in-neighbors, and $(1 \pm
\epsilon) np^2$ common out-neighbors.  Since $\epsilon np^2 \gg 1$,
this implies that $|S_{a,b}| = (1 \pm 3\epsilon) n^2 p^4$.  Let $N_{a,b}$
be the number of these pairs that, after the permutation $\sigma$,
have the additional properties that \textbf{(i)} both $x$ and $y$ are
on the opposite side of the bipartition from $a$, and \textbf{(ii)}
$y$ is the successor of $x$.

Expose the indices $i, j$ for which $v_i = a$ and $v_j = b$; the
remainder is then an $(n-2)$-permutation over the leftover positions.
A similar calculation to Claim 2 shows that $\E{N_{a,b}} = (1 \pm
4\epsilon) \frac{np^4}{2}$, and a similar concentration argument shows
that $N_{a,b}$ is within $\epsilon \cdot \frac{np^4}{2}$ of its mean
with probability $1 - o(n^{-3})$.  Therefore, a union bound over all
$a,b$ completes the proof of this final claim, and the proof of Lemma
\ref{lem:gamma} \hfill $\Box$

\vspace{3mm}

At this point, we could immediately apply the results of Section
\ref{sec:graph} to pack $\Gamma$ with perfect matchings, which then
correspond to Hamilton cycles in $D$.  However, it is unfortunate that
this would miss most of the edges of $D$, since not all edges of $D$
are in correspondence with edges of $\Gamma$.  The solution is to
iterate Procedure 1 several times.  Some care must be taken because
the edges of $D$ will be covered many times by different $\Gamma$.  We
will specify how to deal with this in Procedure 2 below, but first we
collect two intermediate results which control this multiplicity over
several iterations of Procedure 1.


\begin{lemma}
  \label{lem:gamma-cover}
  Suppose $n$ and $\epsilon$ satisfy $\epsilon n \gg 1$.  Let
  $D$ be an $(\epsilon, p)$-uniform digraph on $n$ vertices ($n$
  even), and randomly and independently construct $r = \frac{2 \cdot
    10^5 \log n}{\epsilon^2 p}$ graphs $\Gamma_1, \ldots, \Gamma_r$
  according to Procedure 1.  Let $\DG_i$ be their corresponding
  digraphs.  Then, with probability $1 - o(n^{-1})$, every edge of $D$
  is covered $(1 \pm 1.03\epsilon) \frac{10^5 \log n}{\epsilon^2}$
  times by the $\DG_i$.
\end{lemma}

\noindent \textbf{Proof.}\, Fix an arbitrary edge $\rt{uv}$, and let
$q$ be the probability that it appears in $\DG$ if $\Gamma$ is
constructed according to Procedure 1.  Since the $\DG_i$ are
independent, the number of times $\rt{uv}$ is covered will then be
$\bin{r, q}$, and we will use the Chernoff bounds to prove
concentration.  So, let $\sigma = (v_1, \ldots, v_n)$ be the
permutation which produces $\Gamma$, and let $A$ and $B$ be the two
sides of the bipartition as defined in Procedure 1.

There are two ways that $\rt{uv}$ can appear in $\DG$: either $u \in
A$ and $v \in B$ and $\rt{vu^+} \in D$ for the successor $u^+$ of $u$,
or $u \in B$ and $v \in A$ and $\rt{v^- u} \in D$ for the predecessor
$v^-$ of $v$.  These two cases are clearly disjoint, so $q = q_1 +
q_2$, where $q_1$ and $q_2$ are the respective probabilities in these
two cases.  To calculate $q_1$, the probability that $u \in A$ is
precisely $1/2$.  Conditioned on this, the probability that $v \in B$
is precisely $\frac{n/2}{n-1}$.  By property (i) of $(\epsilon,
p)$-uniformity, $d^+(v) = (1 \pm \epsilon) np$, so the probability
that some out-neighbor of $v$ becomes the successor of $u$ is
$\frac{(1 \pm \epsilon) np}{n-2}$.  Therefore,
\begin{displaymath}
  q_1 
  = 
  \frac{1}{2} \cdot \frac{n/2}{n-1} \cdot \frac{(1 \pm \epsilon) np}{n-2}
  =
  (1 \pm 1.01\epsilon) \frac{p}{4},
\end{displaymath}
since $\epsilon \gg \frac{1}{n}$.
An analogous calculation shows that $q_2 = (1 \pm 1.01\epsilon)
\frac{p}{4}$ as well.  Hence $q = (1 \pm 1.01\epsilon) \frac{p}{2}$.

Now, the number of times that $\rt{uv}$ is covered by the $\DG_i$ is
precisely $\bin{r, q}$.  Since $rq = (1 \pm 1.01\epsilon) \frac{10^5
  \log n}{\epsilon^2}$, the Chernoff bound implies that the
probability that the Binomial deviates from its mean by more than a
factor of $1 \pm 0.01\epsilon$ is at most
\begin{displaymath}
2  \exp\left\{- \frac{(0.01\epsilon)^2}{3} \cdot (1 - 1.01\epsilon) \frac{10^5 \log n}{\epsilon^2} \right\}
  =
  o(n^{-3}).
\end{displaymath}
Taking a union bound over all $O(n^2)$ edges $\rt{uv}$, we obtain the
desired result.  \hfill $\Box$

\vspace{3mm}

\begin{lemma}
  \label{lem:gamma-rare-consec}
  Suppose $r \ll n^{1/2}$.  Construct $r$ independent $\Gamma_i$
  according to Procedure 1.  Then with probability $1 - o(n^{-1})$,
  every pair of distinct vertices $a, b \in D$ has the property that
  $b$ is the successor of $a$ in at most 5 of the permutations for
  $\Gamma_i$.
\end{lemma}

\noindent \textbf{Proof.}\, For a fixed pair $(a, b)$, the probability
that $b$ is the successor of $a$ in a single run of Procedure 1 is
exactly $\frac{1}{n-1}$.  Therefore, the probability of this occurring
more than 5 times in $r$ independent runs is at most $\binom{r}{6}
\frac{1}{(n-1)^{6}} \leq o(n^{-3})$.  Taking a union bound over all
pairs $(a, b)$, we obtain the result.  \hfill $\Box$

\vspace{3mm}

Since a single run of Procedure 1 packs very few edges into Hamilton
cycles, we must repeat the procedure multiple times, deleting the
packed edges from $D$ after each round.  However, it is important to
maintain the pseudo-random properties through the iterations, and the
set of removed edges after a single run of Procedure 1 would be too
sparse to control the changes in the pseudo-random counts.  The
following extension provides one way to achieve this, by repeating
Procedure 1 enough times to uniformly involve all edges of $D$.

\vspace{3mm}

\noindent \textbf{Procedure 2.}\, This takes as input a digraph $D$
with an even number of vertices, and an integer parameter $r$.
\begin{enumerate}

\item Independently generate the random undirected bipartite graphs
  $\Gamma_1$, \ldots, $\Gamma_r$ according to Procedure 1, and let
  $\sigma_i$, $\DG_i$, $A_i$, and $B_i$ be their corresponding
  permutations, digraphs, and bipartitions, respectively.

\item For each edge $\rt{uv} \in D$, let $I_{\rt{uv}}=\{ i :
  \rt{uv}\text{ is covered by }\DG_i \}$.  If $I_{\rt{uv}}\neq
  \emptyset$, then independently select a uniformly random index in
  $I_{\rt{uv}}$ to label $\rt{uv}$ with.

\item For each $\Gamma_i$, define a subgraph $\Gamma_i'$ by keeping
  each edge $ab$ ($a \in A_i$ and $b \in B_i$) if and only if both
  $\rt{ab}$ and $\rt{ba^+}$ are labeled by $i$, where $a^+$ is the
  successor of $a$ according to $\sigma_i$.

\item For each $i$, let $\DG_i'$ be the digraph containing all
  $\{\rt{a b}, \rt{b a^+}\}$, for each $a b \in \Gamma_i'$, with $a
  \in A_i$, where $a^+$ is the successor of $a$ according to
  $\sigma_i$.
\end{enumerate}

Observe that the final step ensures that the $\DG_i'$ are all
disjoint.  So, Lemma \ref{lem:hamilton-matching} shows that we can
work independently on each $\DG_i'$, packing Hamilton cycles by
packing perfect matchings in $\Gamma_i'$ with Proposition
\ref{prop:graph}.  After this, we will remove all $\DG_i'$ from $D$,
and show that their distribution is sufficiently uniform for us to
maintain the necessary pseudo-random properties.  The following Lemma
shows that the $\Gamma_i'$ themselves are sufficiently pseudo-random
for us to apply Proposition \ref{prop:graph}.


\begin{lemma}
  \label{lem:gamma-prime}
  Suppose $n$, $p$, and $\epsilon$ satisfy $\epsilon^{10} n p^8 \gg
  \log^5 n$.  Let $D$ be an $(\epsilon, p)$-uniform digraph on $n$
  vertices ($n$ even), and conduct Procedure 2 with $r = \frac{2 \cdot
    10^5 \log n}{\epsilon^2 p}$.  Let $\kappa = \frac{10^5 \log
    n}{\epsilon^2}$.  Then, with probability $1 - o(n^{-1})$, every
  $\Gamma_i'$ satisfies the following properties:
  \begin{description}
  \item[(i)] All degrees are $(1 \pm 12\epsilon) \frac{n}{2}
    \big(\frac{p}{\kappa}\big)^2$.

  \item[(ii)] All codegrees between pairs of vertices on the same side
    of the bipartition are $(1 \pm 12\epsilon) \frac{n}{2}
    \big(\frac{p}{\kappa}\big)^4$.
  \end{description}
\end{lemma}

\noindent \textbf{Proof.}\, Our restrictions on $\epsilon, n, p$ allow
us to apply Lemmas \ref{lem:gamma} and \ref{lem:gamma-cover}, so we
have that with probability $1 - o(n^{-1})$, after Step 1 of Procedure
2:
\begin{description}
\item[(a)] Every $\Gamma_i$ has all degrees $(1 \pm 5\epsilon)
  \frac{np^2}{2}$ and all same-side codegrees $(1 \pm 5\epsilon)
  \frac{np^4}{2}$.
\item[(b)] Every edge in $D$ is covered $(1 \pm 1.03\epsilon) \kappa$
  times by the $\DG_i$.
\end{description}

Condition on the above outcome of Step 1, and consider an arbitrary
$\Gamma_i'$, which is derived from $\Gamma_i$ with bipartition $A_i
\cup B_i$.  It suffices to show that in the randomness of Step 2, with
probability $1 - o(n^{-2})$ each individual $\Gamma_i'$ has the
desired properties, since $r = o(n)$ by the given restrictions on
$\epsilon, n, p$.  There are four cases to consider: degrees in $A_i$,
degrees in $B_i$, codegrees in $A_i$, and codegrees in $B_i$.
Fortunately, they will all follow by essentially the same argument.

We begin with the degree of an arbitrary vertex $u \in A_i$.  By (a)
above, the degree of $u$ in $\Gamma_i$ is $d_u = (1 \pm 5\epsilon)
\frac{np^2}{2}$.  Note that this actually corresponds to exactly $d_u$
pairs of directed edges in $D$, of the form $\rt{ux}, \rt{x u_+}$, and
all $2 d_u$ directed edges involved are distinct.  Step 2 of Procedure
2 assigns labels to all directed edges, and the degree of $u$ in
$\Gamma_i'$ is precisely the number of the above pairs of directed
edges for which both edges are labeled $i$.  Since every directed edge
is covered $(1 \pm 1.03\epsilon) \kappa$ times by (b) and the edges
are labeled independently, the probability that both of a given pair
above receive label $i$ is $[ (1 \pm 1.03\epsilon) \kappa ]^{-2}$.
Therefore, the expected degree of $u$ in $\Gamma_i'$ is $(1 \pm
8\epsilon) \frac{n}{2} \big(\frac{p}{\kappa}\big)^2$.  Furthermore,
since all directed edges involved are distinct, the Chernoff bound
shows that the probability of the degree deviating from its mean by
more than a factor of $1 \pm \epsilon$ is at most
\begin{displaymath}
2  \exp\left\{
    - \frac{\epsilon^2}{3} \cdot 
    (1 - 8\epsilon) \frac{n}{2} \left( \frac{p}{\kappa} \right)^2
  \right\}
  \leq
  o(n^{-3}),
\end{displaymath}
because $\epsilon^2 n p^2 / \kappa^2 \gg \log n$.  Therefore, with
probability $1 - o(n^{-3})$, the degree of $u$ in $\Gamma_i'$ is $(1
\pm 10\epsilon) \frac{n}{2} \big(\frac{p}{\kappa}\big)^2$.  Taking a
union bound over all $u \in A_i$ establishes part (i) for those
degrees.

The arguments for the other three parts of the lemma are
similar.  For a vertex $v \in B_i$, its incident edges in $\Gamma_i$
correspond to $d_v$ disjoint pairs of directed edges, so the exact
same argument as above produces the bound for these degrees.

The codegree of a given pair of vertices $u, v \in A_i$ corresponds to
disjoint quadruples of distinct directed edges.  From (a), there are
$(1 \pm 5\epsilon) \frac{np^4}{2}$ such quadruples.  The probability
that a given quadruple is completely labeled by $i$ is $[ (1 \pm
1.03\epsilon) \kappa ]^{-4}$ by (b), so the expected codegree in
$\Gamma_i'$ is $(1 \pm 10\epsilon) \frac{n}{2} \big( \frac{p}{\kappa}
\big)^4$.  By the Chernoff bound, the probability that the codegree
deviates from its expectation by a factor of more than $1 \pm
\epsilon$ is at most
\begin{displaymath}
2  \exp\left\{
    - \frac{\epsilon^2}{3} \cdot 
    (1 - 10\epsilon) \frac{n}{2} \left( \frac{p}{\kappa} \right)^4
  \right\}
  \leq
  o(n^{-4}),
\end{displaymath}
since $\epsilon^2 n p^4 / \kappa^4 \gg \log n$.  Taking a union bound
over all pairs of vertices $u, v \in A_i$ produces the desired bound
for their codegrees.  The argument for codegrees in $B_i$ is similar.
\hfill $\Box$

\vspace{3mm}

The next few lemmas build up to a result which controls how $D$ is
affected by the deletion of all edges in the $\DG_i'$.  The first one
controls the ``first-order'' effect of the deletion process.

\begin{lemma}
  \label{lem:d-delete-out-nbrs}
  Condition on the first step of Procedure 2 covering every edge of
  $D$ $(1 \pm 1.03\epsilon) \kappa$ times by the $\DG_i$.  Fix any vertex
  $a \in D$ and any set $E$ of edges incident to $a$ (all oriented in
  the same way with respect to $a$).  Suppose that $\epsilon^2 |E| /
  \kappa^2 \gg \log n$.  Then with probability $1 - o(n^{-5})$, the
  number of edges of $E$ which are covered by the $\DG_i'$ is
  $(1 \pm 1.05\epsilon) \frac{|E|}{\kappa}$,
\end{lemma}

\noindent \textbf{Proof.}\, We start with the case when all edges in
$E$ are directed out of $a$.  Let the random variable $N$ be the
number of these edges which are covered by the $\DG_i'$.  Since we
conditioned on the first step of Procedure 2, the only remaining
randomness is in the independent assignments of the edge labels.  We
expose these labels in three stages.  First, expose the labels of $E$.
Now each edge $e \in E$ has its label $l(e)$, so we can identify its
partner edge in $\DG_{l(e)}$ which must also receive $l(e)$ in order
for both to remain in $\DG_{l(e)}'$; let $F$ be the set of all partner
edges found in this way.  Next, expose all labels outside $E \cup F$.
Finally, expose the labels of $F$.  It is clear that the second stage
does not affect $N$ at all.  Thus, after conditioning on the result of
the second stage, we are at the following situation: every edge $e \in
E$ has a label $l(e)$, and it will only still be covered by the
$\DG_i'$ if its partner in $\DG_{l(e)}$ also receives the label
$l(e)$.  Since its partner is in $(1 \pm 1.03\epsilon) \kappa$
different $\DG_i$, the probability that $e$ is covered by $\DG_i'$ is
the inverse of this multiplicity; linearity of expectation then gives
$\E{N} = (1 \pm 1.04\epsilon) \frac{|E|}{\kappa}$.

We will use the Hoeffding-Azuma inequality to show the concentration
of $N$.  Indeed, the third stage exposure is a product space of
dimension $|F| \leq |E|$.  Consider the effect of changing the label
of a single edge $f \in F$.  By the definition of $F$, the edge $f$ is
either some $\rt{xa}$ directed into $a$, or some $\rt{xb}$ not
incident to $a$ with $\rt{ax} \in E$.  In the latter case, $\rt{ax}$
is the only edge of $E$ which can be affected by the label of $f$, so
$N$ can change by at most 1.  For the remaining case $f = \rt{xa}$,
suppose that the label of $f$ was changed from $j$ to $k$.  By
construction, only the (single) partner edge of $f$ in $\DG_j$ could
suffer from changing $f$'s label away from $j$; this could decrease
$N$ by at most 1.  Similarly, only the partner edge of $f$ in $\DG_k$
could benefit from changing $f$'s label to $k$, and this would only
increase $N$ by at most 1.  We conclude that $N$ is 1-Lipschitz, so
the Hoeffding-Azuma inequality implies that the probability $N$
deviates from its expectation by over $0.01\epsilon
\frac{|E|}{\kappa}$ is at most
\begin{displaymath}
  2 \exp\left\{
    -\frac{ (0.01 \epsilon |E|/\kappa)^2 }{2 |E|}
  \right\}
  \leq
  o(n^{-5}).
\end{displaymath}
Therefore, with probability $1 - o(n^{-5})$, we have $N = (1 \pm
1.05\epsilon) \frac{|E|}{\kappa}$, in the case when all edges of $E$ are
directed out of $a$.  The case when all edges are directed into $a$
follows by essentially the same argument.  \hfill $\Box$

\vspace{3mm}

The next lemma controls the ``second-order'' effect of the deletion
process.

\begin{lemma}
  \label{lem:d-delete-co-out}
  Fix any distinct vertices $a, b \in D$ and any set $X$ of vertices
  such that either
  \begin{description}
  \item[(i)] for all $x \in X$, $\rt{ax}, \rt{bx}$ are edges of $D$; or
  \item[(ii)] for all $x \in X$, $\rt{xa}, \rt{xb}$ are edges of $D$; or
  \item[(iii)] for all $x \in X$, $\rt{ax}, \rt{xb}$ are edges of $D$.
  \end{description}
  Condition on the first step of Procedure 2 covering every edge of
  $D$ at least $\kappa$ times by the $\DG_i$, and on the fact that $b$
  is the successor of $a$ in at most 5 of the permutations.

  Suppose that $|X| / \kappa^4 \gg \log n$.  Then with probability $1
  - o(n^{-5})$, the number $N$ of vertices in $X$ which still have
  both of their designated edges above still covered by the $\DG_i'$
  is at most $\frac{8}{\kappa^2} |X|$.
\end{lemma}

\noindent \textbf{Proof in situations (i) and (ii).}\, These two cases
follow by very similar arguments; for concreteness, let us begin with
(i).  Let $E$ be the set of $2 |X|$ edges specified in (i).  Since we
conditioned on the first step of Procedure 2, the only remaining
randomness is in the independent assignments of the edge labels.  We
expose these labels in three stages.  First, expose the labels of $E$.
Now each edge $e \in E$ has its label $l(e)$, so we can identify its
partner edge in $\DG_{l(e)}$; denote the partner by $\phi(e)$, and let
$F$ be the set of all partner edges found in this way.  Next, expose
all labels outside $E \cup F$.  Finally, expose the labels of $F$.  It
is clear that the second stage does not affect $N$ at all.  Thus,
after conditioning on the result of the second stage, we are at the
following situation: every edge $e \in E$ has a label $l(e)$, and $e$
will only still be covered by the $\DG_i'$ if its partner $\phi(e)$
also receives the label $l(e)$.

The random variable $N$ counts the number of vertices $x$ for which
both $\rt{ax}$ and $\rt{bx}$ are still covered by the $\DG_i'$.  Note
that their partner edges $\phi(\rt{ax})$ and $\phi(\rt{bx})$ are
distinct unless they both coincide as some $\rt{xc}$.  If they are
distinct, then clearly the probability of having both receive their
correct labels is at most $\frac{1}{\kappa^2}$.  On the other hand, if
they coincide, then we must have that in both permutations
$l(\rt{ax})$ and $l(\rt{bx})$, the vertex $c$ is the successor of $a$
and $b$.  Yet $c$ has a unique predecessor, so $l(\rt{ax}) \neq
l(\rt{bx})$.  Therefore, it is actually impossible for the (coincident)
partner edge to receive a label which suits both $\rt{ax}$ and
$\rt{bx}$, so the probability is 0.  In all cases, we have an upper
bound of $\frac{1}{\kappa^2}$, so $\E{N} \leq \frac{|X|}{\kappa^2}$.

We use the Hoeffding-Azuma inequality to probabilistically bound $N$,
since the third exposure stage is a product space.  Consider the
effect of changing the label of a single edge $f \in F$.  By the
definition of $F$, the edge $f$ is either some $\rt{xa}$ or $\rt{xb}$,
or some $\rt{xc}$ with $\rt{ax}$ and $\rt{bx}$ both in $D$.  In the
latter case, $x$ is the only vertex which may have its count in $N$
affected, so $N$ changes by at most 1 under this perturbation.  By
symmetry, it remains to consider the case when $f = \rt{xa}$ changes
its label from $j$ to $k$.  By construction, only the partner edge
$\rt{ay}$ of $f$ in $\DG_j$ could suffer from changing $f$'s label
away from $j$; this could potentially lose only $y$ in the count of
$N$, so $N$ would decrease by at most 1.  Similarly, only the partner
edge of $f$ in $\DG_k$ could benefit from changing $f$'s label to $k$,
and this could only increase $N$ by at most 1.  We conclude that $N$
is 1-Lipschitz over its product space of dimension $|F| \leq 2 |X|$,
so the Azuma-Hoeffding inequality implies that the probability that
$N$ exceeds its expectation by more than $\frac{|X|}{\kappa^2}$ is at
most
\begin{displaymath}
2  \exp\left\{
    - \frac{ ( |X| / \kappa^2 )^2}{ 4|X| }
  \right\}
  \leq
  o(n^{-5}).
\end{displaymath}
This finishes the case when (i) holds.  The case when (ii) holds
follows from a directly analogous argument.  

\vspace{3mm}

\noindent \textbf{Proof in situation (iii).}\, Let $E$ be the set of
$2 |X|$ edges specified in (iii).  We expose the labels in the same
three stages as before: first $E$, then the other non-partner edges,
and finally the partner edges $F$.  

Again, after the second stage we would be at the following situation:
every edge $e \in E$ has a label $l(e)$, and a vertex $x$ will only be
counted toward $N$ if the partner edge $\phi(\rt{ax})$ receives the
label $l(\rt{ax})$, and $\phi(\rt{xb})$ also receives $l(\rt{xb})$.
However, this time there is an additional complication, because it can
happen that one or both of these partner edges have already had their
labels exposed.  The only way this could happen is if either
$\phi(\rt{ax}) = \rt{xb}$ or $\phi(\rt{xb}) = \rt{ax}$.

Fortunately, we are only seeking an upper bound on $N$, so we only
need to determine when this first round exposure already forces a
vertex $x$ to contribute to $N$.  A moment's thought reveals that the
only way this can happen is if the first round gave both $\rt{ax}$ and
$\rt{xb}$ the same label $l$, and furthermore, the permutation
$\sigma_l$ has $b$ as the successor of $a$.  Now we use the assumption
that at most 5 of the permutations $\sigma_l$ satisfy this property.
Let $L$ be the set of the corresponding indices $l$.

We can circumvent this issue, by observing that in the first exposure,
the number of $x$ for which both $\rt{ax}$ and $\rt{xb}$ receive the
same label $l\in L$ is stochastically dominated by $\bin{|X|,
  5/\kappa^2}$.  Since $|X| / \kappa^2 \gg \log n$, the Chernoff bound
then implies that with probability $1 - o(n^{-5})$, the first round
has at most $\frac{6}{\kappa^2} |X|$ ``bad'' $x$ which have the potential
of being automatically included in the count for $N$.

The second round exposure is essentially irrelevant, so we may now
condition on the result of the second round satisfying the property in
the previous paragraph.  As we only need to upper bound $N$, it
remains to consider only the ``good'' $x$.  From our previous
discussion, if a good $x$ has, say, the label of $\phi(\rt{ax})$
already exposed to be $l$ then $b$ is not the successor of $a$ in
permutation $l$ and so this $x$ cannot contribute to $N$.

Also note that the only way for $\phi(\rt{ax})$ to equal
$\phi(\rt{xb})$ is if both are $\rt{ba}$.  If $\phi(\rt{ax})=\rt{ba}$
then $a\in B$ and if $\phi(\rt{xb})=\rt{ba}$ then $a\in A$ and so we
must have $l(\phi(\rt{ax})) \neq l(\phi(\rt{xb}))$.  Hence it is not
possible for $\rt{ba}$ to simultaneously match both labels, and $x$
cannot contribute to $N$.

So, if we let let $Y$ be the subset of good vertices $x \in X$ for
which $\phi(\rt{ax}) \neq \phi(\rt{xb})$, and neither label has yet
been exposed, it remains to control the number $M$ of vertices in $Y$
which contribute to $N$.  Since each $x \in Y$ has $\phi(\rt{ax}) \neq
\phi(\rt{xb})$ and their labels are independent, we immediately have
$\E{M} \leq |Y| / \kappa^2$.

We show that $M$ is concentrated by using the Hoeffding-Azuma
inequality on the third round exposure product space of dimension $|F|
\leq 2|X|$.  Consider the effect of changing the label of an edge $f
\in F$ from $j$ to $k$.  There are three types of edges in $F$:
\begin{description}
\item[Case 1: $f = \rt{za}$.] These $f$ can only arise as partner
  edges of some $\rt{ax}$ (possibly several).  However, $f$ has a
  unique partner edge $\rt{ax}$ in $\DG_j$, so changing its label away
  from $j$ can only reduce $M$ at $x$.  Similarly, changing it to $k$
  can only grow $M$ at $y$, where $\rt{ay}$ is the unique partner edge
  of $f$ in $\DG_k$.  Therefore, $M$ can only change by at most 1.

\item[Case 2: $f = \rt{bz}$.] By the same argument as Case 1, this can
  only change $M$ by at most 1.

\item[Case 3: $f = \rt{xy}$.] These $f$ can only arise as partners of
  either $\rt{ax}$ or $\rt{yb}$.  Clearly, changing $l(f)$ can only
  affect whether $x$ or $y$ are counted in $M$, so it has an effect of
  at most 2.
\end{description}
Thus $M$ is 2-Lipschitz, and the Hoeffding-Azuma inequality shows that
the probability it exceeds its expectation by $|X| / \kappa^2 $ is at
most
\begin{displaymath}
  2\exp\left\{
    - \frac{ ( \kappa^{-2} |X| )^2 }{ 4 |X| }
  \right\}
  \leq
  o(n^{-5}).
\end{displaymath}

Since we conditioned on there being at most $\frac{6}{\kappa^2} |X|$
bad vertices, we have $N \leq M + \frac{6}{\kappa^2} |X|$.  Putting
everything together, we conclude that with probability $1 -
o(n^{-5})$, $N \leq \frac{8}{\kappa^2} |X|$, completing the proof.
\hfill $\Box$

\vspace{3mm}

The previous two lemmas now enable us to prove that the pseudo-random
properties of the digraph are still maintained after deleting the
$\DG_i'$.

\begin{lemma}
  \label{lem:digraph-delete}
  Suppose $n$, $p$, and $\epsilon$ satisfy $\epsilon^8 np^4 \gg \log^5
  n$.  Let $r = \frac{2 \cdot 10^5 \log n}{\epsilon^2 p}$ and $\kappa
  = \frac{10^5 \log n}{\epsilon^2}$.  Let $D$ be an $(\epsilon,
  p)$-uniform digraph on $n$ vertices ($n$ even), and conduct
  Procedure 2.  Let $D'$ be the subgraph of $D$ obtained by deleting
  all edges covered by any $\DG_i'$.  Then, with probability $1 -
  o(n^{-1})$, $D'$ is an $(\epsilon', p')$-uniform digraph with
  $\epsilon' = \epsilon \big( 1 + \frac{4.23}{\kappa} \big)$ and $p' =
  p \big( 1 - \frac{1}{\kappa} \big)$.
\end{lemma}

\noindent \textbf{Proof.}\, By applying Lemmas \ref{lem:gamma-cover}
and \ref{lem:gamma-rare-consec}, the outcome of Step 1 of Procedure 2
will satisfy the following two properties with probability $1 -
o(n^{-1})$.  (The second follows since $(\epsilon^2 p)^4 \gg n^{-1}$,
and so $r \ll n^{1/2}$.)
\begin{itemize}
\item Every edge of $D$ covered $(1 \pm 1.03\epsilon) \kappa$ times by
  the $\DG_i$.
\item For any $a,b$, the vertex $b$ is the successor of $a$ in at most
  5 of the permutations for the $\Gamma_i$.
\end{itemize}
Condition on this outcome.  We will now show that in the remaining
randomness of Step 2, the properties of $(\epsilon', p')$-uniformity
are satisfied with probability $1 - o(n^{-1})$.  We establish them one
at a time.

\vspace{3mm}

\noindent \textbf{Property (i).}\, Consider an arbitrary vertex $a$.  
By $(\epsilon, p)$-uniformity, it has $(1 \pm \epsilon) np$
out-edges in $D$, and $\epsilon^2 (np) / \kappa^2 = \Theta(\epsilon^6
np / \log^2 n) \gg \log n$.  So, we may apply Lemma
\ref{lem:d-delete-out-nbrs} to the set of out-edges of $a$.  This
shows that with probability $1 - o(n^{-5})$, the new out-degree of $a$
in $D'$, in terms of its original out-degree $d^+(a)$ in $D$, is
\begin{eqnarray*}
  d^+(a) - (1 \pm 1.05\epsilon) \frac{d^+(a)}{\kappa}
  &=&
  d^+(a) \left[
    1 - (1 \pm 1.05\epsilon) \frac{1}{\kappa}
  \right] \\
  &=&
  (1 \pm \epsilon) np \cdot \left(
    1 \pm \frac{1.06\epsilon}{\kappa}
  \right)
  \left(
    1 - \frac{1}{\kappa}
  \right) \\
  &=& (1 \pm \epsilon') np'.
\end{eqnarray*}
Here, we used $d^+(a) = (1 \pm \epsilon) np$ by $(\epsilon,
p)$-uniformity.  Taking a union bound over all $a \in D$, we obtain
the desired result on out-degrees with probability $1 - o(n^{-4})$.  A
similar argument controls all in-degrees.  \hfill $\Box$

\vspace{3mm}

\noindent \textbf{Property (ii).}\, Consider an arbitrary pair of
distinct vertices $a, b$, and let $X$ be the set of their common
out-neighbors.  Let $X_1$ be the number of vertices in $x \in X$ such
that $\rt{ax}$ is covered by the $\DG_i'$, let $X_2$ be the number of
$x \in X$ such that $\rt{bx}$ is covered by the $\DG_i'$, and let
$X_{12}$ be the number of $x \in X$ such that both $\rt{ax}$ and
$\rt{bx}$ are covered by the $\DG_i'$.  Clearly, the number of common
out-neighbors of $a$ and $b$ in $D'$ is exactly $d^+(a, b) - X_1 - X_2
+ X_{12}$, where $d^+(a, b)$ was the number of their common out-neighbors
in $D$.

Note that $d^+(a, b) = (1 \pm \epsilon) np^2$ by $(\epsilon,
p)$-regularity, and $\epsilon^2 (np^2) / \kappa^2 = \Theta(\epsilon^6
np^2 / \log^2 n) \gg \log n$.  So, Lemma \ref{lem:d-delete-out-nbrs}
implies that with probability $1 - o(n^{-5})$, both $X_1$ and $X_2$
are $(1 \pm 1.05\epsilon) \frac{ d^+(a, b) }{\kappa}$.  On the other
hand, we also have $(np^2) / \kappa^4 = \Theta( \epsilon^8 np^2 /
\log^4 n ) \gg \log n$, so Lemma \ref{lem:d-delete-co-out} bounds
$X_{12}$ by $\frac{9}{\kappa^2} d^+(a, b)$ with probability $1 -
o(n^{-5})$.  Therefore, $X_{12}$ is within $0.01\epsilon \frac{d^+(a,
  b)}{\kappa}$ additive error of $\frac{1}{\kappa^2} d^+(a, b)$
because $\epsilon \kappa \gg 1$.

Putting these bounds together, we have that the new number of common
out-neighbors is
\begin{eqnarray*}
  d^+(a, b) - X_1 - X_2 + X_{12} &=&
  d^+(a, b) - 2 \cdot (1 \pm 1.05\epsilon) \frac{1}{\kappa} d^+(a, b) 
  + \left(
    \frac{1}{\kappa^2} d^+(a, b) \pm \frac{0.01\epsilon}{\kappa} d^+(a, b)
  \right) \\
  &=&
  d^+(a, b) \left[
    1 - \frac{2}{\kappa} + \frac{1}{\kappa^2} \pm \frac{2.11 \epsilon}{\kappa}
  \right] \\
  &=&
  (1 \pm \epsilon) np^2 
  \cdot
  \left(
    1 \pm \frac{2.12\epsilon}{\kappa}
  \right)
  \left(
    1 - \frac{1}{\kappa}
  \right)^2 \\
  &=&
  (1 \pm \epsilon') n(p')^2,
\end{eqnarray*}
Here, $d^+(a, b) = (1 \pm \epsilon) np^2$ by $(\epsilon,
p)$-uniformity.  Taking a union bound over all $a, b \in D$, we obtain
the desired result on the new $d^+(a, b)$ with probability $1 -
o(n^{-3})$.  Similar arguments control the other two parts of property
(ii) of $(\epsilon', p')$-uniformity.  \hfill $\Box$

\vspace{3mm}

\noindent \textbf{Property (iii).}\, This is a slight extension of the
previous argument.  Consider any four vertices $a, b, c, d$, which are
all distinct except for the possibility $b=c$.  Let $X$ and $X'$ be
the sets of vertices $x$ such that $\rt{ax}$, $\rt{xb}$, $\rt{cx}$,
$\rt{xd}$ are all in $D$ or $D'$, respectively.  Let $X_1$, $X_2$,
$X_3$, and $X_4$ be the sets of vertices $x \in X$ such that
$\rt{ax}$, $\rt{xb}$, $\rt{cx}$, or $\rt{xd}$ are still covered by the
$\DG_k'$.  Define the pairwise intersections $X_{ij} = X_i \cap X_j$.
By inclusion-exclusion, we have
\begin{displaymath}
  |X| - |X_1| - |X_2| - |X_3| - |X_4| 
  \ \ \leq \ \ 
  |X'|
  \ \ \leq \ \ 
  |X| - |X_1| - |X_2| - |X_3| - |X_4| 
  + \sum_{i < j} |X_{ij}|.
\end{displaymath}
As in (ii), since $\epsilon^2 (np^4) / \kappa^2 = \Theta(\epsilon^6
np^4 / \log^2 n) \gg \log n$ (to apply Lemma
\ref{lem:d-delete-out-nbrs}) and $(np^4) / \kappa^4 =
\Theta(\epsilon^8 np^4 / \log^4 n) \gg \log n$ (for Lemma
\ref{lem:d-delete-co-out}), we have that with probability $1 -
o(n^{-5})$, each $|X_i| = (1 \pm 1.05\epsilon) \frac{1}{\kappa} |X|$,
and each $|X_{ij}| \leq \frac{9}{\kappa^2} |X| = o\big(
\frac{\epsilon}{\kappa}|X| \big)$.  Thus
\begin{displaymath}
  |X'|
  =
  |X| \left[
    1 - \frac{4}{\kappa} 
    \pm 4 \cdot \frac{1.05\epsilon}{\kappa}
    \pm \frac{0.01 \epsilon}{\kappa}
  \right]
  = 
  (1 \pm \epsilon) np^4 \left( 
    1 - \frac{4}{\kappa} \pm \frac{4.21\epsilon}{\kappa} 
  \right).
\end{displaymath}
Here, we used $|X| = (1 \pm \epsilon) np^4$ by $(\epsilon,
p)$-uniformity.  Now observe that
\begin{displaymath}
  \left(
    1 - \frac{1}{\kappa}
  \right)^4
  =
  1 - \frac{4}{\kappa} + O\left( \frac{1}{\kappa^2} \right).
\end{displaymath}
However, since $\epsilon \kappa \gg 1$, the error term is
$o\big( \frac{\epsilon}{\kappa} \big)$.  Therefore, we have
\begin{displaymath}
  |X| 
  = 
  (1 \pm \epsilon) np^4 \left( 
    1 - \frac{4}{\kappa} \pm \frac{4.21\epsilon}{\kappa} 
  \right)
  =
  (1 \pm \epsilon) np^4 \left( 
    1 \pm \frac{4.22\epsilon}{\kappa} 
  \right)
  \left(
    1 - \frac{1}{\kappa}
  \right)^4
  =
  (1 \pm \epsilon') n (p')^4.
\end{displaymath}
Taking a union bound over all choices of $a, b, c, d$, we obtain the
desired result on the new $d^+(a, b)$ with probability $1 -
o(n^{-1})$.  \hfill $\Box$

\vspace{3mm}

Now we have established control over the pseudo-random properties of
the digraph after deletion, as well as over the individual bipartite
graphs $\Gamma_i'$ at each stage.  We finally combine all of our
lemmas to prove the main result of this section, that pseudo-random
digraphs can be efficiently packed with Hamilton cycles.

\vspace{3mm}

\noindent \textbf{Proof of Theorem \ref{thm:digraph}.}\, We will
iterate Procedure 2, packing all intermediate $\DG_i'$ at each
iteration and deleting them from the digraph, until very few edges
remain.  The total number of iterations will be $o(n)$, so since our
previous lemmas hold with probability $1 - o(n^{-1})$, a final union
bound will show that we achieve the desired packing \whp.

Let $D_0 = D$, $\epsilon_0 = \epsilon$, and $p_0 = p$.  Define the
sequences $(\epsilon_t)$, $(p_t)$ via the following recursion:
\begin{eqnarray*}
  \epsilon_{t+1} &=& \epsilon_t \left( 1 + \frac{4.23 \epsilon_t^2}{10^5 \log n} \right) \\
  p_{t+1} &=& p_t \left( 1 - \frac{\epsilon_t^2}{10^5 \log n}
  \right).
\end{eqnarray*}
Let $T$ be the smallest index such that $p_T \leq \frac{1}{8}
\epsilon^{1/8} p$.  Note that $\epsilon_t$ only increases, so
\begin{displaymath}
  T 
  <
  \frac{10^5 \log n}{\epsilon^2} \cdot  \log \frac{1}{\epsilon}
  \ll
  n.
\end{displaymath}
Also note that since 
\begin{displaymath}
  (1 + 4.23x) (1 - x)^{4.23} 
  \leq 
  e^{4.23x} \big( e^{-x} \big)^{4.23}
  = 
  1,
\end{displaymath}
we have in general that $\frac{\epsilon_{t+1}}{\epsilon_t} \leq \big(
\frac{p_t}{p_{t+1}} \big)^{4.23}$.  Therefore, 
\begin{equation}
  \label{eq:epsilon-T}
  \epsilon_{T-1} 
  \leq 
  \epsilon \cdot \big( 8 \epsilon^{-1/8} \big)^{4.23} 
  \ll 
  \epsilon^{0.47}, 
  \andsp
  \text{and so}
  \andsp
  \epsilon_{T-1}^{1/3} \ll \epsilon^{0.15} \ll \epsilon^{1/8}.
\end{equation}

We now iteratively construct $D_1, \ldots, D_T$, such that each $D_t$
is $(\epsilon_t, p_t)$-uniform.  Indeed, consider the (random)
Procedure 2 applied to $D_t$ with respect to $r_t = \frac{2 \cdot 10^5
  \log n}{\epsilon_t^2 p_t}$.  Let $\kappa_t = \frac{10^5 \log
  n}{\epsilon_t^2}$.  This produces bipartite graphs $\Gamma_{t,i}'$
and directed graphs $\DG_{t,i}'$ with all $\DG_{t,i}'$ disjoint.  Let
$D_{t+1}$ be the result of deleting all $\DG_{t,i}'$ from $D_t$.  To
apply Lemmas \ref{lem:gamma-prime} and \ref{lem:digraph-delete}, we
must check that $\epsilon_t^{10} n p_t^8 \gg \log^5 n$.  But this
follows from our initial assumption that $\epsilon^{11} n p^8 \gg
\log^5 n$ since $\epsilon_t \geq \epsilon$ and $p_t \geq \frac{1}{8}
\epsilon^{1/8} p$.  Therefore, as the individual failure probabilities
of Lemmas \ref{lem:gamma-prime} and \ref{lem:digraph-delete} are
$o(n^{-1})$ and $T = o(n)$, we may assume that the outcome of
Procedure 2 satisfies the following properties:
\begin{description}
\item[(i)] All degrees in each $\Gamma_{t,i}'$ are $(1 \pm 12\epsilon_t)
  \frac{n}{2} \big(\frac{p_t}{\kappa_t}\big)^2$.

\item[(ii)] All codegrees between pairs of vertices on the same side
  of the bipartition of each $\Gamma_{t,i}'$ are $(1 \pm 12\epsilon_t)
  \frac{n}{2} \big(\frac{p_t}{\kappa_t}\big)^4$.

\item[(iii)] $D_{t+1}$ is $(\epsilon_{t+1}, p_{t+1})$-uniform.
\end{description}
Fix this outcome.  We may apply Proposition \ref{prop:graph} to each
$\Gamma_{t,i}'$ because
\begin{displaymath}
  \epsilon_t^{4/3} n \left( \frac{p_t}{\kappa_t} \right)^4 
  >
  \Theta
  \left(
    \epsilon_t^{4/3} n \cdot (\epsilon^{1/8} p)^4 \cdot \frac{\epsilon_t^8}{\log^4 n}
  \right)
  \gg
  \epsilon^{10} n p^4 / \log^4 n
  \gg 1,
\end{displaymath}
where we used $p_t \geq \frac{1}{8} \epsilon^{1/8} p$, $\epsilon_t
\geq \epsilon$, and our initial assumption that $\epsilon^{11} np^8
\gg \log^5 n$.  So, every $\Gamma_{t,i}'$ can be packed with
edge-disjoint perfect matchings, missing only a
$4\epsilon_t^{1/3}$-fraction of the edges.  By Lemma
\ref{lem:hamilton-matching}, these edge-disjoint perfect matchings in
$\Gamma_{t,i}'$ correspond to edge-disjoint Hamilton cycles in
$\DG_{t,i}'$, missing the same fraction of edges since there is a
2-to-1 correspondence between edges in $\DG_{t,i}'$ and
$\Gamma_{t,i}'$.

We carry on the above procedure until we create $D_T$.  Then, we will
have packed Hamilton cycles in $D \setminus D_T$, up to a fractional
error of $4 \epsilon_{T-1}^{1/3}$.  Since $D_T$ is $(\epsilon_T,
p_T)$-uniform, it has at most $(1 + \epsilon_T) n^2 p_T \leq 2 n^2
\cdot \frac{1}{8} \epsilon^{1/8} p$ edges.  As $D$ was $(\epsilon,
p)$-uniform, it had at least $(1 - \epsilon) n^2 p \geq \frac{1}{2}
n^2 p$ edges.  Thus $D_T$ itself has at most $\frac{1}{2}
\epsilon^{1/8}$-fraction of the total number of edges.

Therefore, the fraction of edges of $D$ that have not
been covered is at most
\begin{displaymath}
  4 \epsilon_{T-1}^{1/3} \cdot \left( 1 - \frac{1}{2} \epsilon^{1/8} \right)
  +
  \frac{1}{2} \epsilon^{1/8}
  <
  \epsilon^{1/8},
\end{displaymath}
since $\epsilon_{T-1}^{1/3} \ll \epsilon^{1/8}$ by inequality
\eqref{eq:epsilon-T}.
This completes the proof.

\vspace{3mm}

\noindent \textbf{Remark.}\, The above proof showed that we can pack
the edges of $D$ with Hamilton cycles, up to a fractional error of at
most $\epsilon^c$, where $c = 1/8$.  In fact, one can prove the same
result for any $c > 1/7$, and all of the calculations can be recovered
by following our proof.  Indeed, we intentionally introduced a very
high level of precision in Lemma \ref{lem:gamma-cover}.  It is clear
that the same argument can replace the $1.03$ with $1.003$, etc.
These smaller order errors are what accumulate into the decimal places
in the final $4.23$ which appears in Lemma \ref{lem:digraph-delete},
and which later determines the value of $c$.  When $4.23$ is replaced
by a constant arbitrarily close to 4, the above argument produces a
result with $c$ approaching $1/7$.

\section{Hypergraphs}
\label{sec:3graph}

This section is organized so that its structure is parallel to the
previous section.  Indeed, their themes are similar, as the objective
of this section is to establish another reduction, this time from the
hypergraph setting to the digraph setting.  The fundamental building
block is now the following procedure, which should be compared with
Procedure 1.

\vspace{3mm}

\noindent \textbf{Procedure 3.}\, This takes as input a 3-graph $H$
with an even number of vertices.
\begin{enumerate}
\item Generate a random permutation $\sigma = (v_1, v_2, \ldots, v_n)$
  of the vertices of $H$.  Split it into $n/2$ successive pairs $(v_1,
  v_2)$, $(v_3, v_4)$, \ldots, $(v_{n-1}, v_n)$.

\item Construct a directed graph $D$ with $n/2$ vertices, one
  corresponding to each of the pairs above.  Place a directed edge
  from $(v_i, v_{i+1})$ to $(v_j, v_{j+1})$ if and only if both
  hyperedges $e=\{v_i, v_{i+1}, v_j\}$ and $f=\{v_{i+1}, v_j,
  v_{j+1}\}$ are present in $H$.  Consider $e$ and $f$ to be
  \emph{partners}, and introduce the notation $\f_\sigma(e)=f$ and
  $\f_\sigma(f)=e$.

\item Since a directed edge in $D$ corresponds to two hyperedges in
  $H$, we account for this by defining $H' \subset H$ to be the
  hypergraph containing all $\{v_i, v_{i+1}, v_j\}$, $\{v_{i+1}, v_j,
  v_{j+1}\}$, for each directed edge from $(v_i, v_{i+1})$ to $(v_j,
  v_{j+1})$ in $D$.
\end{enumerate}

\vspace{3mm}

As already noted, Hamilton cycles in $D$ correspond precisely to tight
Hamilton cycles in $H$.  Indeed, the hyperedges corresponding to the
directed edges of a Hamilton cycle are precisely those which are
needed for a tight Hamilton cycle in the hypergraph.

\begin{lemma}
  \label{lem:hamilton-di-hyp}
  Let $C_1, C_2$ be edge-disjoint Hamilton cycles in $D$.  Then their
  associated Hamilton cycles $C_1', C_2'$ in $H$ are also
  edge-disjoint.
\end{lemma}

\noindent \textbf{Proof.}\, It suffices to show that given a fixed
pairing $(v_1, v_2)$, \ldots, $(v_{n-1}, v_n)$ from Step 1, any given
hyperedge $\{a, b, c\}$ can be associated with at most one directed
edge in Step 3.  Indeed, by construction, some pair of vertices in the
hyperedge, say $(a, b)$, must be a pair above.  If the pair containing
$c$ is $(c, x)$, then the directed edge must be from $(a, b)$ to $(c,
x)$, hence unique.  Similarly, if the pair containing $c$ is $(x, c)$,
then the only possibility is that the directed edge is from $(x, c)$
to $(a, b)$.  \hfill $\Box$

\vspace{3mm}

As in the previous section, our next goal is to show that the result
of Procedure 3 (a digraph) is pseudo-random.  The following lemma will
be a useful component of the proof.

\begin{lemma}
  \label{lem:pairs-kept}
  Let $S$ be a given set of ordered pairs of distinct vertices $(a, b)$,
  with $\epsilon^2 |S|^2 / n^3 \gg \log n$.  Then with probability $1
  - o(n^{-9})$, the number $N$ of ordered pairs in $S$ which appear
  in the list of $n/2$ ordered pairs from Step 1 of Procedure 3,
  satisfies $N = (1 \pm \epsilon) \frac{|S|}{2n}$.
\end{lemma}

\noindent \textbf{Proof.}\, The probability that any given $(x,y) \in
S$ appears as an ordered pair in Step 1 is precisely $\frac{1}{2}
\cdot \frac{1}{n-1}$.  This is because $x$ must land in one of the
$n/2$ positions of the permutation which correspond to first-vertices
of ordered pairs, and then $y$ must land in the position immediately
following $x$.  Since $\epsilon \gg \frac{\sqrt{n^3 \log n}}{|S|} \gg
\frac{1}{\sqrt{n}}$, linearity of expectation gives $\E{N} = \big( 1
\pm \frac{\epsilon}{2} \big) \frac{|S|}{2n}$.

For concentration, we apply Fact \ref{fact:conc-perm}.  Suppose the
permutation $\sigma$ is converted to $\sigma'$ via the single
transposition of the positions of the vertices $x$ and $y$.  Let $x'$
and $y'$ be the other vertices in their respective ordered pairs
(according to $\sigma$).  Then, the only changes in the set of ordered
pairs produced by Step 1 is that the ordered pairs $(x, x')$ and $(y,
y')$ change to $(y, x')$ and $(x, y')$.  All other ordered pairs
remain the same, so their inclusion or exclusion in $N$ remains
invariant.  Therefore, the value of $N$ can only change by at most 2.
Fact \ref{fact:conc-perm} then implies that the probability that $N$
deviates from its mean by more than $\frac{\epsilon}{2}
\frac{|S|}{2n}$ is at most
\begin{displaymath}
2  \exp\left\{  
    -\frac{ 2 \left(
        \frac{\epsilon}{2} \frac{|S|}{2n}
      \right)^2 }
    {2^2 n}
  \right\}
  \leq
  o(n^{-9}),
\end{displaymath}
completing the proof.  \hfill $\Box$.

\vspace{3mm}

Now we show that the result of Procedure 3 is a pseudo-random digraph.

\begin{lemma}
  \label{lem:d}
  Suppose $n$, $p$, and $\epsilon$ satisfy $\epsilon^2 n p^{16} \gg
  \log n$.  Let $H$ be an $(\epsilon, p)$-uniform 3-graph on $n$
  vertices ($n$ even), and randomly construct the $\frac{n}{2}$-vertex
  directed graph $D$ according to Procedure 3.  Then, with probability
  $1-o(n^{-1})$, $D$ is $(5\epsilon, p^2)$-uniform.
\end{lemma}

\begin{figure}[htbp]
  \centering
  \begin{tabular}{c @{\hspace{1cm}} c @{\hspace{1cm}} c @{\hspace{1cm}} c @{\hspace{1cm}} c}
    \includegraphics[scale=2]{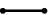} & 
    \includegraphics[scale=2]{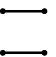} & 
    \includegraphics[scale=2]{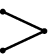} & 
    \includegraphics[scale=2]{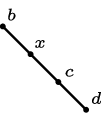} & 
    \includegraphics[scale=2]{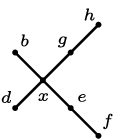} \\
    $\Gamma_1$ & $\Gamma_2$ & $\Gamma_3$ & $\Gamma_4$ & $\Gamma_5$ \\
  \end{tabular}
  \caption{Auxiliary graphs}
  \label{fig:gamma}
\end{figure}

\noindent \textbf{Proof.}\, We verify the properties of $D$ one at a
time, starting with out-degrees.  Fix any two distinct vertices $a,
b$, and let $N_{a,b}$ be the number of ordered pairs $(x, y)$ produced
by Step 1 such that $\{a, b, x\}, \{b, x, y\} \in H$.  It suffices to
show that with probability $1 - o(n^{-3})$, $N_{a,b} = (1 \pm
5\epsilon) \frac{n}{2} p^2$.  (Note the factor of 2 because $D$ has
$n/2$ vertices.)  For this, first let $S_{a,b}$ be the set of ordered
pairs $(x, y)$, not necessarily produced by Step 1, such that $\{a, b,
x\}, \{b, x, y\} \in H$.  This is a deterministic set, from the
hypergraph $H$, whose size we can estimate by applying $(\epsilon,
p)$-uniformity to $H$ with various auxiliary graphs $\Gamma$.  Indeed,
by applying the property with $\Gamma_1$ (in Figure \ref{fig:gamma})
and $(a, b)$, the number of choices for $x$ such that $\{a, b, x\} \in
H$ is $(1 \pm \epsilon) np$.  Then, given such a choice for $x$,
applying uniformity with $\Gamma_1$ and $(b, x)$ shows that there are
$(1 \pm \epsilon) np$ choices for $y$ such that $\{b, x, y\}$ is also
an edge.  Putting these estimates together, we conclude that
$|S_{a,b}| = (1 \pm 3\epsilon) n^2 p^2$.  Since $\epsilon^2 (n^2
p^2)^2 / n^3 \gg \log n$, Lemma \ref{lem:pairs-kept} gives $N_{a,b} =
(1 \pm 5\epsilon) \frac{n}{2} p^2$ with the desired probability.

The arguments for the other properties are very similar.  For
in-degrees, we fix a pair of distinct vertices $(a, b)$, and apply
uniformity with $\Gamma_1$ to show there are $(1 \pm \epsilon) np$
choices for $y$ such that $\{y, a, b\} \in H$.  Applying again on $(y,
a)$, we find a further $(1 \pm \epsilon) np$ choices for $x$ such that
$\{x, y, a\} \in H$.  This is exactly the same estimate as we had in
the out-degree case above, so the same argument completes it.

To control the common out-degrees of pairs of vertices in $D$, we fix 4 distinct vertices of $H$,
arranged in two pairs $(a, b)$ and $(c, d)$.  Applying uniformity with
$\Gamma_2$ (see Figure \ref{fig:gamma}), there are $(1 \pm \epsilon)
np^2$ choices for $x$ such that $\{a, b, x\}, \{c, d, x\} \in H$.
After fixing such a choice for $x$, another application of uniformity
with $\Gamma_3$ (see Figure \ref{fig:gamma}) shows that there are $(1
\pm \epsilon) np^2$ choices for $y$ such that $\{b, x, y\}, \{d, x,
y\} \in H$.  Therefore, there are $(1 \pm 3\epsilon) n^2 p^4$ pairs
$(x, y)$ which could serve as the terminal vertex of both directed
edges from $(a, b)$ and $(c, d)$, if only $(x, y)$ appeared as a pair
in Step 1.  Since $\epsilon^2 (n^2 p^4)^2 / n^3 \gg \log n$, we may
apply Lemma \ref{lem:pairs-kept}, and finish as before.

A similar argument controls $d^-(a, b)$ in $D$.  To establish the last
part of property (ii) of $(5\epsilon, p^2)$-regularity for $D$, fix
any 4 distinct vertices in $H$ arranged in two pairs $(a, b)$ and $(c,
d)$.  Applying uniformity with $\Gamma_1$ (see Figure
\ref{fig:gamma}), there are $(1 \pm \epsilon) np$ choices for $x$ such
that $\{a, b, x\} \in H$.  Fix such an $x$.  Applying uniformity again
with $\Gamma_4$ (see Figure \ref{fig:gamma}), with $(b, x, c, d)$
corresponding to the labels in Figure \ref{fig:gamma}, we see that
there are $(1 \pm \epsilon) np^3$ choices for $y$ such that all of
$\{b, x, y\}, \{x, y, c\}, \{y, c, d\} \in H$.  Hence there are $(1
\pm 3\epsilon) n^2 p^4$ choices for a pair $(x, y)$ which would be
both an in-neighbor of $(a, b) \in D$ and an out-neighbor of $(c, d)
\in D$, if only $(x, y)$ appeared as a pair in Step 1.  Thus Lemma
\ref{lem:pairs-kept} finishes this case as before.

For property (iii) of $(5\epsilon, p^2)$-regularity for $D$, fix any 8
distinct vertices in $H$, arranged in 4 pairs $(a, b)$, $(c, d)$, $(e,
f)$, and $(g, h)$.  Applying uniformity with $\Gamma_2$ shows that
there are $(1 \pm \epsilon) np^2$ choices for $x$ such that both $\{a,
b, x\}, \{c, d, x\} \in H$.  Then, after fixing such an $x$, applying
uniformity with $\Gamma_5$ (see Figure \ref{fig:gamma}) shows that
there are $(1 \pm \epsilon) np^6$ choices for $y$ such that all $\{b,
x, y\}, \{x, y, e\}, \{y, e, f\}, \{d, x, y\}, \{x, y, g\}, \{y, g,
h\} \in H$.  Hence there are $(1 \pm 3\epsilon) n^2 p^8$ choices for a
pair $(x, y)$ which would out-neighbors of both $(a, b), (c, d) \in D$
and in-neighbors of $(e, f), (g, h) \in D$, if only $(x, y)$ appeared
as a pair in Step 1.  Since $\epsilon^2 (n^2 p^8)^2 / n^3 \gg \log n$,
we may apply Lemma \ref{lem:pairs-kept}, and finish as before.  \hfill
$\Box$

\vspace{3mm}

Next we prove the analogue of Lemma \ref{lem:gamma-cover},
establishing that repeated runs of Procedure 3 will cover the edges of
$H$ fairly uniformly.

\begin{lemma}
  \label{lem:d-cover}
  Suppose $n$ and $\epsilon$ satisfy $\epsilon n \gg 1$.  Let
  $H$ be an $(\epsilon, p)$-uniform 3-graph on $n$ vertices ($n$
  even), and randomly and independently construct $r = \frac{ 10^6 n
    \log n}{3 \epsilon^2 p}$ digraphs $D_1, \ldots, D_r$ according to
  Procedure 3.  Let $H_i$ be their corresponding 3-graphs.  Then, with
  probability $1 - o(n^{-1})$, every edge of $H$ is covered $(1 \pm
  1.03\epsilon) \frac{10^6 \log n}{\epsilon^2}$ times by the $H_i$.
\end{lemma}

\noindent \textbf{Proof.}\, Fix an edge $\{a, b, c\}$ of $H$, and let
$q$ be the probability that it is covered by $H_1$.  There are several
ways in which this can happen.  First, it could be that $(a, b)$ is a
pair from the permutation $\sigma_1$ for $H_1$, and $(c, x)$ is
another pair, with $\{b, c, x\}$ also an edge of $H$.  Let us bound
the probability $q_1$ of this particular event happening.  Let $t$ be
the number of choices for $x$ such that $\{b, c, x\} \in H$.  In terms
of $t$, we have $q_1 = \frac{1}{2} \cdot \frac{1}{n-1} \cdot
\frac{1}{2} \cdot \frac{t}{n-3}$.  To see this, first expose the
position of $a$; it is in the first half of an ordered pair with
probability $1/2$.  Next, $b$ must be occupy the second position in
that ordered pair, and then $c$ must select a first half spot in one
of the remaining ordered pairs.  Finally, there are $t$ valid choices
for $c$'s partner.  But by $(\epsilon, p)$-uniformity applied with
$\Gamma_1$ (Figure \ref{fig:gamma}), $t = (1 \pm \epsilon) np$.  So
(using $\epsilon \gg \frac{1}{n}$), $q_1 = (1 \pm 1.01\epsilon)
\frac{p}{4n}$.

The above argument can be run with any permutation of $a, b, c$.
Also, it could be that $c$ occupied the second position in its ordered
pair, e.g., if $(x,c)$ and $(a,b)$ were the pairs involved.
Therefore, $q = 2 \cdot 3! \cdot (1 \pm 1.01\epsilon) \frac{p}{4n} =
(1 \pm 1.01\epsilon) \frac{3p}{n}$.  Since the $r$ random
constructions are independent, the number of times $\{a, b, c\}$ is
covered is $\bin{r, q}$, which has expectation $rq = (1 \pm
1.01\epsilon) \frac{10^6 \log n}{\epsilon^2}$.  The Chernoff bound
implies that the probability that the Binomial deviates from its mean
by more than a factor of $1 \pm 0.01\epsilon$ is at most
\begin{displaymath}
2  \exp\left\{- \frac{(0.01\epsilon)^2}{3} \cdot (1 - 1.01\epsilon) \frac{10^6 \log n}{\epsilon^2} \right\}
  =
  o(n^{-4}).
\end{displaymath}
Taking a union bound over all $O(n^3)$ edges $\{a, b, c\}$, we obtain
the desired result.  \hfill $\Box$

\vspace{3mm}

Next, we prove a result which will serve the same purpose as Lemma
\ref{lem:gamma-rare-consec} of the previous section.  Procedure 3
constructs a digraph $D$ on half as many vertices by pairing up the
vertices of $H$.  Let us say that a set $S$ of 4 vertices of $H$ is
\emph{condensed}\/ in $D$ if two of those pairs above contain all 4
vertices in $S$.  For example, if $S = \{a, b, c, d\}$, and two
vertices of $D$ are the ordered pairs $(b, c)$ and $(d, a)$, then $S$
is condensed in $D$.

\begin{lemma}
  \label{lem:d-rare-close}
  Suppose $r \ll n^{3/2}$.  Construct $r$ independent $D_i$ according
  to Procedure 3.  Then with probability $1 - o(n^{-1})$, every set
  $S$ of 4 vertices in $H$ is condensed in at most 9 of the $D_i$.
\end{lemma}

\noindent \textbf{Proof.}\, For a fixed set of 4 vertices $\{a, b, c,
d\}$, the probability that it is condensed in $D_1$ is precisely
$\frac{3}{n-1} \cdot \frac{1}{n-3} < \frac{4}{n^2}$.  To see this,
first expose the position of $a$ in the permutation $\sigma_1$ which
defines $D_1$.  This identifies the position which contains the other
vertex in $a$'s pair; expose the vertex in that spot.  The probability
that it is $b$, $c$, or $d$ is precisely $\frac{3}{n-1}$.  Assuming
this is the case, there are still two vertices of interest left, say
$c$ and $d$.  Then expose the position of $c$.  This identifies the
position which contains the other vertex in $c$'s pair; expose the
vertex in that spot.  The probability that it is $d$ is precisely
$\frac{1}{n-3}$.

Since the $D_i$ are independent, the number of them which have the
above property with respect to $\{a, b, c, d\}$ is stochastically
dominated by $\bin{r, \frac{4}{n^2}}$.  Since we assumed $r \ll
n^{3/2}$, the probability that this exceeds 9 is at most
\begin{displaymath}
  \binom{r}{10} \left( \frac{4}{n^2} \right)^{10}
  =
  o(n^{-5}),
\end{displaymath}
so a union bound over all $O(n^4)$ choices for $\{a, b, c, d\}$
completes the proof.  \hfill $\Box$

\vspace{3mm}

We now formulate the analogue of Procedure 2 for the hypergraph
setting.  

\vspace{3mm}

\noindent \textbf{Procedure 4.}\, This takes as input a 3-graph $H$
with an even number of vertices, and an integer parameter $r$.
\begin{enumerate}

\item Independently generate the random digraphs $D_1, \ldots, D_r$
  according to Procedure 3, and let $H_i$ be their corresponding
  3-graphs.

\item For each edge $e \in H$ let $I_e = \{ i : e\text{ is covered by
  }H_i \}$.  If $I_e \neq \emptyset$, independently select a uniformly
  random index in $I_e$ to label $e$ with.

\item For each $D_i$, define the subgraph $D_i'$ by keeping each edge
  $\rt{uv}$ if and only if both of its corresponding 3-graph edges are
  labeled by $i$.

\item For each $i$, let $H_i'$ be the 3-graph containing all
  hyperedges which correspond to the directed edges left in $D_i'$.
  Note that $e(H_i') = 2 e(D_i')$.
\end{enumerate}

As in Procedure 2, the final step ensures that the $H_i'$ are all
disjoint, so Lemma \ref{lem:hamilton-di-hyp} shows that we can work on
each $D_i'$ independently, packing Hamilton cycles with Theorem
\ref{thm:digraph}.  Then, we will remove all $H_i'$ from $H$, and
repeat.  The next result establishes the necessary pseudo-random
properties of the $D_i'$.

\begin{lemma}
  \label{lem:d-prime}
  Suppose $n$, $p$, and $\epsilon$ satisfy $\epsilon^{18} n p^{16} \gg
  \log^9 n$.
  Let $H$ be an $(\epsilon, p)$-uniform 3-graph on $n$ vertices ($n$
  even), and conduct Procedure 4 with $r = \frac{10^6 n \log n}{3
    \epsilon^2 p}$.  Let $\kappa = \frac{10^6 \log n}{\epsilon^2}$.
  Then, with probability $1 - o(n^{-1})$, every $D_i'$ is $\big(
  16\epsilon, \frac{p^2}{\kappa^2} \big)$-uniform.
\end{lemma}

\noindent \textbf{Proof.}\, Our restrictions on $\epsilon, n, p$ allow
us to apply Lemma \ref{lem:d} ($\epsilon^2 n p^{16} \gg \log n$) and
Lemma \ref{lem:d-cover} $(\epsilon \gg n^{-1})$, so we have that with
probability $1 - o(n^{-1})$, after Step 1 of Procedure 4:
\begin{description}
\item[(a)] Every $D_i$ is $(5\epsilon, p^2)$-uniform.
\item[(b)] Every edge in $H$ is covered $(1 \pm 1.03\epsilon) \kappa$
  times by the $H_i$.
\end{description}

Condition on the above outcome of Step 1, and consider an arbitrary
$D_i'$.  It suffices to show that in the randomness of Step 2, with
probability $1 - o(n^{-3})$ each individual $D_i'$ has the desired
properties, since $r = o(n^2)$ by the given restrictions on $\epsilon,
n, p$.  We start by verifying out-degrees.  A vertex $v \in D_i'$
corresponds to a pair of vertices $(a, b)$ from $H$.  For each edge of
$D_i$ which is directed away from $v$, there is a distinct pair of
vertices $(x, y)$ from $H$.  That directed edge remains in $D_i'$ if
and only if both of its associated hyperedges $\{a, b, x\}, \{b, x,
y\}$ receive label $i$.  This happens with probability $[(1 \pm
1.03\epsilon) \kappa]^{-2}$ by property (b) above.  Yet property (a)
ensures that the number $t$ of such directed edges in $D_i$ is $(1 \pm
5\epsilon) \frac{n}{2} p^2$, so the expected out-degree of $v$ in
$D_i'$ is $(1 \pm 8\epsilon) \frac{n}{2} \frac{p^2}{\kappa^2}$.

For concentration, note that all of the $2t$ hyperedges associated
with the $t$ directed edges are distinct, so their labels are
generated independently.  Therefore, by the Chernoff inequality, the
probability that the out-degree of $v$ in $D_i'$ deviates from its
expectation by more than a factor of $1 \pm \epsilon$ is at most
\begin{displaymath}
2  \exp\left\{
    - \frac{\epsilon^2}{3} \cdot 
    (1 - 8\epsilon) \frac{n}{2} \frac{p^2}{\kappa^2}
  \right\}
  \leq
  o(n^{-4}),
\end{displaymath}
because $\epsilon^2 n p^2 / \kappa^2 = \Theta(\epsilon^6 n p^2 /
\log^2 n) \gg \log n$.  Therefore, with probability $1 - o(n^{-4})$,
the degree of $v$ in $D_i'$ is $(1 \pm 10\epsilon) \frac{n}{2}
\frac{p^2}{\kappa^2}$.  Taking a union bound over all $v \in D_i$
establishes uniformity for out-degrees.

All other properties of $\big( 16\epsilon, \frac{p^2}{\kappa^2}
\big)$-uniformity follow by a similar argument.  Importantly, in each
case, all directed edges correspond to disjoint pairs of hyperedges,
so we can always use the Chernoff inequality to establish
concentration.  The smallest mean we ever deal with is $(1 \pm
5\epsilon) \frac{n}{2} \big( \frac{p^2}{(1 \pm 1.02\epsilon)\kappa^2}
\big)^4 = (1 \pm 10\epsilon) \frac{n}{2} \big( \frac{p^2}{\kappa^2}
\big)^4$, from property (iii) of uniformity for $D_i'$, so the largest
error in concentration is at most
\begin{displaymath}
2  \exp\left\{
    - \frac{\epsilon^2}{3} \cdot 
    (1 - 10\epsilon) \frac{n}{2} \left( \frac{p^2}{\kappa^2} \right)^4
  \right\}
  \leq
  o(n^{-7}),
\end{displaymath}
where we use $\epsilon^2 n (p^2 / \kappa^2)^4 = \epsilon^{18} n p^8 /
\log^8 n \gg \log n$.  \hfill $\Box$

\vspace{3mm}

We now move to show how the hypergraph $H$ is affected by the deletion
of all hypergraphs $H_i'$.  As before, we use inclusion/exclusion to
sandwich the quantities in question, using an accurate ``first-order''
estimate, together with a rough ``second-order'' upper bound.

\begin{lemma}
  \label{lem:h-delete-1}
  Condition on the first step of Procedure 4 covering every edge of
  $H$ $(1 \pm 1.03\epsilon) \kappa$ times by the $H_i$.  Fix any two
  distinct vertices $a, b \in H$, and a set $S$ of vertices such that
  $\{a, b, x\} \in H$ for all $x \in S$.  Suppose that $\epsilon^2 |S|
  / \kappa^3 \gg \log n$.  Then with probability $1 - o(n^{-8})$, the
  number $N$ of vertices $x \in S$ such that $\{a, b, x\}$ is covered
  by the $H_i'$ is $(1 \pm 1.05\epsilon) \frac{|S|}{\kappa}$,
\end{lemma}

\noindent \textbf{Proof.}\, Let $S = \{x_1, \ldots, x_t\}$.  Since we
are conditioning on the result of Step 1, the only randomness left is
in the independent exposure of all hyperedge labels, most of which are
irrelevant for $N$.  We define $F$, the set of relevant edges, as
follows.  For each $H_j$ which covers some $e=\{a, b, x_i\}$, there is
a partner edge $f=\phi_j(e)$ such that $e$ remains in $H_j'$ if and
only if $e$ and $f$ receive the label $j$.  Let $F$ be the collection
of all partner edges obtained in this way, together with all $\{a, b,
x_i\}$.  Since we assume every $\{a, b, x_i\}$ is covered at most $(1
+ 1.03\epsilon) \kappa$ times by the $H_j$, we always have $|F| \leq 2
\kappa |S|$.  Observe that all labels outside $F$ have no effect on
$N$ at all, so we may condition on an arbitrary setting of those
labels, leaving only $|F|$ independent labels to be exposed.

Next, we show that the probability that a particular edge $e=\{a, b,
x_i\}$ is covered by the $H_j'$ is $[ (1 \pm 1.03\epsilon) \kappa
]^{-1}$.  To see this, note that after revealing the label $j$ of $e$,
there is a single partner edge $f=\phi_j(e)$, whose label must match
that of $e$.  Since $f$ is covered $(1 \pm 1.03\epsilon) \kappa$ times
by the $H_i$, it has probability $[ (1 \pm 1.03\epsilon) \kappa
]^{-1}$ of receiving that label $j$.  Therefore by linearity of
expectation, $\E{N} = (1 \pm 1.04\epsilon) \frac{|S|}{\kappa}$.

We use the Hoeffding-Azuma inequality to establish the concentration
of $N$, as the remaining randomness is a product space of dimension
$|F| \leq 2 \kappa |S|$.  Consider the effect of changing the label of
a hyperedge $e = \{x, y, z\} \in F$ from $i$ to $j$.  Note that its
old partner edge $\phi_i(e)$ overlapped with $e$ in two vertices, but
also included a new vertex, say $w$.  Changing the label of $e$ away
from $i$ could potentially remove any or all of $x, y, z, w$ from the
count for $N$ (if $\phi_i(e)$ had received label $i$), so this could
potentially reduce $N$ by up to 4.  Similarly, setting the label of
$e$ to $j$ can potentially increase $N$ by up to 4.  Therefore, $N$ is
4-Lipschitz, so the Hoeffding-Azuma inequality implies that the
probability that $N$ deviates from its expectation by more than
$0.01\epsilon \frac{|S|}{\kappa}$ is at most
\begin{displaymath}
  2 \exp\left\{
    -\frac{ (0.01 \epsilon |S|/\kappa)^2 }{2 \cdot 4^2 |F|}
  \right\}
  \leq
  2 \exp\left\{
    -\frac{ 0.01^2 \epsilon^2 |S| }{64 \kappa^3}
  \right\}
  \leq
  o(n^{-8}),
\end{displaymath}
where we used $\epsilon^2 |S| / \kappa^3 \gg \log n$.  Therefore, with
probability $1 - o(n^{-8})$, we have $N = (1 \pm 1.05\epsilon)
\frac{|S|}{\kappa}$, as desired.  \hfill $\Box$

\vspace{3mm}

Having proved an accurate ``first-order'' estimate, we move to the
``second-order'' upper bound (analogous to Lemma
\ref{lem:d-delete-co-out}).  In the 3-graph setting, this breaks into
two cases, which we treat in the following two lemmas.

\begin{lemma}
  \label{lem:h-delete-2-disjoint}
  Condition on the first step of Procedure 4 covering every edge of
  $H$ $(1 \pm 1.03\epsilon)\kappa$ times by the $H_i$.  Fix any
  distinct vertices $a, b, c, d \in H$ and any set of vertices $S$
  such that both $\{a, b, x\}, \{c, d, x\} \in H$ for all $x \in S$.
  Suppose that $|S| / \kappa^5 \gg \log n$.  Then with probability $1
  - o(n^{-8})$, the number $N$ of vertices $x \in S$ which still have
  both $\{a, b, x\}$ and $\{c, d, x\}$ covered by the $H_i'$ is at
  most $\frac{3}{\kappa^2} |S|$.
\end{lemma}

\noindent \textbf{Proof.}\, Again, the only remaining randomness is in
the labeling of the hyperedges in Step 2.  Also, since each edge is
covered by at most $(1 + 1.03\epsilon) \kappa$ of the $H_j$, less than
$2 \kappa \cdot 2|S|$ of them are relevant for $N$, for the same
reason as in the previous proof.  So, we may condition on all labels
for irrelevant edges, and focus on the remaining product space over
the relevant edges $F$.

Let us bound the probability that a particular $x \in S$ remains in
the count for $N$.  First, expose the labels for $e = \{a, b, x\}$ and
$f = \{c, d, x\}$, and suppose they are $k$ and $l$.  Let their
partner edges be $e' = \phi_k(e)$ and $f' = \phi_l(f)$.  Importantly,
$e'$ cannot equal $f$, because $e'$ must overlap with $e$ in two
vertices.  Hence the label of $e'$ has not yet been exposed.
Similarly, the label of $f'$ has not yet been exposed.  However, note
that it may happen that $e' = f'$.

In order for $x$ to remain in the count for $N$, $e'$ must receive
label $k$ and $f'$ must receive label $l$.  First we show that if $e'
= f'$, then the probability of this occurring is 0.  Indeed, if $k
\neq l$, it clearly cannot happen.  But otherwise, if $k = l$, it is
impossible for $e' = f'$ in the first place, because in $H_k$, each
edge has a distinct partner edge, so $e'$ and $f'$ should be distinct.

On the other hand, if $e' \neq f'$, the probability that these
independent labels equal $k$ and $l$ is at most $[(1 - 1.03\epsilon)
\kappa ]^{-2}$, because each edge is covered by at least $(1 -
1.03\epsilon) \kappa$ of the $H_j$.  Therefore, we conclude that in
all cases, the probability of a certain $x_i$ remaining in the count
for $N$ is at most $\frac{2}{\kappa^2}$, so $\E{N} \leq
\frac{2}{\kappa^2} |S|$.

Finally, the same argument that we used in the proof of Lemma
\ref{lem:h-delete-1} shows that $N$ is 4-Lipschitz over the product
space of dimension $|F| < 4 \kappa |S|$.  Hence the probability that
$N$ exceeds its expectation by more than $\frac{|S|}{\kappa^2}$ is at
most
\begin{displaymath}
2  \exp\left\{
    - \frac{ ( |S|/\kappa^2 )^2}{ 2 \cdot 4^2 |F| }
  \right\}
  \leq
2  \exp\left\{
    - \frac{|S|}{ 128 \kappa^5 }
  \right\}
  \leq
  o(n^{-8}),
\end{displaymath}
since we assumed $|S| / \kappa^5 \gg \log n$.  Therefore, with
probability $1 - o(n^{-8})$, we have $N \leq \frac{3}{\kappa^2} |S|$,
as desired.  \hfill $\Box$

\vspace{3mm}

\begin{lemma}
  \label{lem:h-delete-2-incident}
  Condition on the first step of Procedure 4 covering every edge
  of $H$ $(1 \pm 1.03\epsilon)\kappa$ times by the $H_i$, and on the
  property that every set of 4 vertices of $H$ is condensed in at most
  9 of the $D_i$.  Fix any distinct vertices $a, b, c \in H$ and any
  set of vertices $S$ such that both $\{a, b, x\}, \{b, c, x\} \in H$
  for all $x \in S$.  Suppose that $|S| / \kappa^5 \gg \log n$.  Then
  with probability $1 - o(n^{-8})$, the number $N$ of vertices $x \in
  S$ which still have both $\{a, b, x\}$ and $\{b, c, x\}$ covered by
  the $H_i'$ is at most $\frac{13}{\kappa^2} |S|$.
\end{lemma}

\noindent \textbf{Proof.}\, Again, the only remaining randomness is in
the labeling of the hyperedges in Step 2, so we may concentrate on the
product space of dimension $|F| \leq 2 \kappa \cdot 2|S|$, where $F$
is the set of all edges $\{a, b, x\}, \{b, c, x\}$ and their possible
partners.

To bound the probability that a particular $x \in S$ remains in the
count for $N$, expose the labels for $e = \{a, b, x\}$ and $f = \{b,
c, x\}$, and suppose they are $k$ and $l$.  Let their partner edges be
$e' = \phi_k(e)$ and $f' = \phi_l(f)$.  This time, it is possible for
$e' = f$ or $f' = e$, but that can only happen if the 4 vertices $a,
b, x, c$ are condensed in $D_k$.  Furthermore, if $k \neq l$, then $x$
is automatically excluded from the count for $N$, which is not a
problem because we only seek an upper bound.

Therefore, the probability that $x$ contributes to $N$ is at most $q_1
+ q_2$, where $q_1$ is the probability that $e$ and $f$ both receive
the same label $k$, with $a, b, x, c$ condensed in $D_k$, and $q_2$ is
the conditional probability that $e'$ and $f'$ receive labels $k$ and
$l$, given the exposure of the labels of $e$ and $f$, and $e' \neq f$,
$f' \neq e$.  We assumed that $a, b, x, c$ were only condensed in at
most 9 of the $D_j$, so $q_1 \leq \frac{9}{(1 - 1.03\epsilon)\kappa}
\cdot \frac{1}{(1 - 1.03\epsilon)\kappa} < \frac{10}{\kappa^2}$.

As for $q_2$, the same argument that we used in the previous lemma
shows that if $e' = f'$, then we must have had $k \neq l$, so it is
impossible for $e'$ to receive $k$ and $f'$ to receive $l$, hence $x$
cannot contribute to $N$.  The only remaining case is $e' \neq f'$,
and the probability that both receive the correct label is at most
$[(1 - 1.03\epsilon) \kappa ]^{-2}$.  Putting everything together, we
conclude that the probability that $x$ contributes to $N$ is at most
$\frac{12}{\kappa^2}$, so $\E{N} \leq \frac{12}{\kappa^2} |S|$.

Finally, the same argument that we used in the previous proof again
establishes the Lipschitz constant of 4 for $N$ in the product space
of dimension $|F| < 4 \kappa |S|$.  Hence the probability that $N$
exceeds its expectation by more than $\frac{|S|}{\kappa^2}$ is at most
\begin{displaymath}
2  \exp\left\{
    - \frac{ ( |S|/\kappa^2 )^2}{ 2 \cdot 4^2 |F| }
  \right\}
  \leq
2  \exp\left\{
    - \frac{|S|}{ 128 \kappa^5 }
  \right\}
  \leq
  o(n^{-8}),
\end{displaymath}
since we assumed $|S| / \kappa^5 \gg \log n$.  Therefore, with
probability $1 - o(n^{-8})$, we have $N \leq \frac{13}{\kappa^2} |S|$,
as desired.  \hfill $\Box$

\vspace{3mm}

Now we combine our first-order and second-order bounds, and show that
the pseudo-random properties of $H$ are maintained after deleting all
$H_i'$.

\begin{lemma}
  \label{lem:h-delete}
  Suppose $n$, $p$, and $\epsilon$ satisfy $\epsilon^{10} np^6 \gg
  \log^6 n$.  Let $r = \frac{10^6 n \log n}{3 \epsilon^2 p}$ and
  $\kappa = \frac{10^6 \log n}{\epsilon^2}$.  Let $H$ be an
  $(\epsilon, p)$-uniform 3-graph on $n$ vertices ($n$ even), and
  conduct Procedure 4.  Let $H'$ be the subgraph of $H$ obtained by
  deleting all edges covered by any $H_i'$.  Then, with probability $1
  - o(n^{-1})$, $H'$ is an $(\epsilon', p')$-uniform 3-graph with
  $\epsilon' = \epsilon \big( 1 + \frac{6.6}{\kappa} \big)$ and $p' =
  p \big( 1 - \frac{1}{\kappa} \big)$.
\end{lemma}

\noindent \textbf{Proof.}\, By applying Lemmas \ref{lem:d-cover} and
\ref{lem:d-rare-close}, the outcome of Step 1 of Procedure 4 will
satisfy the following two properties with probability $1 - o(n^{-1})$.
(The second follows since $(\epsilon^2 p)^5 \gg n^{-1}$, and so $r \ll
n^{3/2}$.)
\begin{itemize}
\item Every edge of $H$ is covered $(1 \pm 1.03\epsilon) \kappa$ times
  by the $H_i$.
\item Every set of 4 vertices of $H$ is condensed in at most 9 of the
  $D_i$.
\end{itemize}
Condition on this outcome.  We will now show that in the remaining
randomness of Step 2, the properties of $(\epsilon', p')$-uniformity
are satisfied with probability $1 - o(n^{-1})$.  For this, fix a set
of $t \leq 7$ vertices $v_1, \ldots v_t \in H$, and let $\Gamma$ be an
arbitrary $t$-vertex graph with $s \leq 6$ edges.  It suffices to show
that in $H'$, the random variable $X = d_\Gamma(v_1, \ldots, v_t)$ is
$(1 \pm \epsilon') np^s$ with probability $1 - o(n^{-8})$, because we
can then take a union bound over all $O(n^7)$ choices for the $v_i$,
and all $O(1)$ possibilities for $\Gamma$.

Let $S$ be the set of all $x$ such that $\{v_i, v_j, x\} \in H$ for
every edge $ij \in \Gamma$.  Our assumed $(\epsilon, p)$-uniformity
gives $|S| = (1 \pm \epsilon) np^s$.  We will use a similar
inclusion/exclusion argument as in the proof of Lemma
\ref{lem:digraph-delete} to estimate $X$ in terms of $|S|$.  Let $e_1,
\ldots, e_s$ be the edges of $\Gamma$.  For each $k \in \{1, \ldots,
s\}$, define the random variable $X_k$ as follows.  Let $i, j$ be the
endpoints of edge $e_k$.  Then let $X_k$ be the number of vertices $x
\in S$ such that $\{v_i, v_j, x\}$ is covered by some $H_l'$.  Also,
for every two distinct $k_1, k_2 \in \{1, \ldots, s\}$, let the random
variable $X_{k_1 k_2}$ count the number of vertices $x \in S$ such
that both $\{v_{i_1}, v_{j_1}, x\}$ and $\{v_{i_2}, v_{j_2}, x\}$ are
covered by some $H_l'$, where $i_1, j_1$ and $i_2, j_2$ are the
respective endpoints of $e_1$ and $e_2$.  In terms of these random
variables, the principle of inclusion/exclusion always gives
\begin{displaymath}
  |S| - \sum_{i=1}^s X_i
  \leq
  X
  \leq
  |S| - \sum_{i=1}^s X_i + \sum_{i<j} X_{ij}.
\end{displaymath}
Since we noted above that $|S| = (1 \pm \epsilon) np^s =
\Omega(np^6)$, we have $\epsilon^2 |S| / \kappa^3 \geq
\Omega(\epsilon^8 np^6 / \log^3 n) \gg \log n$, so Lemma
\ref{lem:h-delete-1} controls all $X_i = (1 \pm 1.05\epsilon)
\frac{|S|}{\kappa}$ with probability $1 - o(n^{-8})$.  Also, since
$|S| / \kappa^5 \geq \Omega(\epsilon^{10} np^6 / \log^5 n) \gg \log
n$, Lemmas \ref{lem:h-delete-2-disjoint} and
\ref{lem:h-delete-2-incident} show that every $X_{ij} <
\frac{13}{\kappa^2} |S|$ with probability $1 - o(n^{-8})$.

Putting everything together, and using $s \leq 6$ and
$\frac{1}{\kappa} \ll \epsilon$, we have
\begin{displaymath}
  X 
  =
  |S| - s(1 \pm 1.05\epsilon) \frac{|S|}{\kappa}
  \pm s^2 \frac{13}{\kappa^2} |S|
  =
  |S| \cdot \left(1 - \frac{s}{\kappa} \pm \frac{6.3\epsilon}{\kappa} \pm \frac{0.1 \epsilon}{\kappa} \right)
  =
  (1 \pm \epsilon) np^s \cdot \left(1 - \frac{s}{\kappa} \pm \frac{6.4\epsilon}{\kappa} \right)
\end{displaymath}
Now observe that
\begin{displaymath}
  \left(
    1 - \frac{1}{\kappa}
  \right)^s
  =
  1 - \frac{s}{\kappa} + O\left( \frac{1}{\kappa^2} \right).
\end{displaymath}
Again using $\frac{1}{\kappa} \ll \epsilon$, we bound the error term
by $o\big( \frac{\epsilon}{\kappa} \big)$.  Therefore, we have
\begin{displaymath}
  X 
  = 
  (1 \pm \epsilon) np^s \left( 
    1 - \frac{s}{\kappa} \pm \frac{6.4\epsilon}{\kappa} 
  \right)
  =
  (1 \pm \epsilon) np^s \left[ 
    \left(
      1 - \frac{1}{\kappa}
    \right)^s
    \pm \frac{6.5\epsilon}{\kappa} 
  \right]
  =
  (1 \pm \epsilon') n (p')^s,
\end{displaymath}
as desired.  \hfill $\Box$

\vspace{3mm}

We finally finish the proof of Theorem \ref{thm:3graph}.  The method
is essentially the same as that used to prove Theorem
\ref{thm:digraph} in the previous section.

\vspace{3mm}

\noindent \textbf{Proof of Theorem \ref{thm:3graph}.}\, Let $H_0
= H$, $\epsilon_0 = \epsilon$, and $p_0 = p$.  Define the sequences
$(\epsilon_t)$, $(p_t)$ via the following recursion:
\begin{eqnarray*}
  \epsilon_{t+1} &=& \epsilon_t \left( 1 + \frac{6.6 \epsilon_t^2}{10^6 \log n} \right) \\
  p_{t+1} &=& p_t \left( 1 - \frac{\epsilon_t^2}{10^6 \log n}
  \right).
\end{eqnarray*}
Let $T$ be the smallest index such that $p_T \leq \frac{1}{2}
\epsilon^{1/15} p$.  Note that $\epsilon_t$ only increases, so
$$  T 
  \leq 
  \frac{10^6 \log n}{\epsilon^2} \cdot  \log \frac{1}{\epsilon}.
$$
Also note that since 
\begin{displaymath}
  (1 + 6.6x) (1 - x)^{6.6} 
  \leq 
  e^{6.6x} \big( e^{-x} \big)^{6.6}
  = 
  1,
\end{displaymath}
we have in general that $\frac{\epsilon_{t+1}}{\epsilon_t} \leq \big(
\frac{p_t}{p_{t+1}} \big)^{6.6}$.  Therefore, 
\begin{equation}
  \label{eq:3epsilon-T}
  \epsilon_{T-1} 
  \leq 
  \epsilon \cdot \big( 2 \epsilon^{-1/15} \big)^{6.6} 
  =
  \Theta(\epsilon^{0.56}), 
  \andsp
  \text{and so}
  \andsp
  \epsilon_{T-1}^{1/8} = \Theta( \epsilon^{0.07} ) \ll \epsilon^{1/15}.
\end{equation}

We now iteratively construct $H_1, \ldots, H_T$, such that each $H_t$
is $(\epsilon_t, p_t)$-uniform.  Indeed, consider the (random)
Procedure 4 applied to $H_t$ with respect to $r_t = \frac{10^6 n \log
  n}{3 \epsilon_t^2 p_t}$.  Let $\kappa_t = \frac{10^6 \log
  n}{\epsilon_t^2}$.  This produces digraphs $D_{t,i}'$ and 3-graphs
$H_{t,i}'$ with all $H_{t,i}'$ disjoint.  Let $H_{t+1}$ be the result
of deleting all $H_{t,i}'$ from $H_t$.  To apply Lemmas
\ref{lem:d-prime} and \ref{lem:h-delete}, we must check that
$\epsilon_t^{18} n p_t^{16} \gg \log^9 n$.  But this follows from our
initial assumption that $\epsilon^{45} np^{16} \gg \log^{21} n$ since
$\epsilon_t \geq \epsilon$ and $p_t \geq \frac{1}{2} \epsilon^{1/15}
p$.  Therefore, the two lemmas show that with probability $1 -
o(n^{-1})$, Procedure 4 results in the following properties:
\begin{description}
\item[(i)] Every $D_{t,i}'$ is $\big( 16 \epsilon_t,
  \frac{p_t^2}{\kappa_t^2} \big)$-uniform.

\item[(ii)] $H_{t+1}$ is $(\epsilon_{t+1}, p_{t+1})$-uniform.
\end{description}
We may assume this outcome, since we will only iterate $T = o(n)$
times.  In order to apply Theorem \ref{thm:digraph} to each
$D_{t,i}'$, we must verify that $\epsilon_t^{11} \frac{n}{2}
\big(\frac{p_t^2}{\kappa_t^2}\big)^8 \gg \log^5 n$.  Indeed, this is
the case:
\begin{displaymath}
  \epsilon_t^{11} \frac{n}{2} \left( \frac{p_t^2}{\kappa_t^2} \right)^8
  \gg
  \Theta \left(
    \epsilon_t^{11} n \cdot (\epsilon^{1/15} p)^{16} 
    \cdot \frac{\epsilon_t^{32}}{\log^{16} n}
  \right)
  \gg
  \epsilon^{45} n p^{16} / \log^{16} n
  \gg
  \log^5 n,
\end{displaymath}
by our initial assumption that $\epsilon^{45} np^{16} \gg \log^{21}
n$.  So, every $D_{t,i}'$ can be packed with Hamilton cycles, missing
only an $\epsilon_t^{1/8}$-fraction of the edges.  By Lemma
\ref{lem:hamilton-di-hyp}, these edge-disjoint Hamilton cycles in
$D_{t,i}'$ correspond to edge-disjoint Hamilton cycles in $H_{t,i}'$,
missing the same fraction of edges since there is a 2-to-1
correspondence between edges in $H_{t,i}'$ and $D_{t,i}'$.

We carry on the above procedure until we create $H_T$.  Then, we will
have packed Hamilton cycles in $H \setminus H_T$, up to error of
$\epsilon_{T-1}^{1/8}$-fraction.  It remains to estimate the numbers
of edges in $H_T$ and $H$.  Note that in general, by applying
$(\epsilon, p)$-uniformity to every pair of vertices (using $\Gamma_1$
in Figure \ref{fig:gamma}), we can estimate the number of edges in any
$(\epsilon, p)$-uniform 3-graph to be $\frac{1}{3!} \cdot n^2 \cdot (1
\pm \epsilon) np$, because this counts every hyperedge $3!$ times.
Thus $H$ has at least $(1 - \epsilon) \frac{n^3 p}{6} \geq \frac{n^3
  p}{7}$ edges, and $H_T$ has at most $(1 + \epsilon_T) \frac{n^3
  p_T}{6} \leq \frac{n^3 \epsilon^{1/15} p}{11}$ edges.  Thus $H_T$
itself has at most $\frac{7}{11} \epsilon^{1/15}$-fraction of the
total number of edges.

Therefore, the fraction of edges of $H$ that have not
been covered is at most
\begin{displaymath}
  \epsilon_{T-1}^{1/8} \cdot \left( 1 - \frac{7}{11} \epsilon^{1/15} \right)
  +
  \frac{7}{11} \epsilon^{1/15}
  \leq
  \epsilon^{1/15},
\end{displaymath}
since $\epsilon_{T-1}^{1/8} \ll \epsilon^{1/15}$ by inequality
\eqref{eq:3epsilon-T}.  This completes the proof.  \hfill $\Box$

\vspace{3mm}

\noindent \textbf{Remark.}\, The above proof showed that we can pack
the edges of $H$ with Hamilton cycles, up to a fractional error of at
most $\epsilon^c$, where $c = 1/15$.  In fact, one can prove the same
result for any $c > 1/13$, and all of the calculations can be
recovered by following our proof.  Just as in Section
\ref{sec:digraph}, we intentionally introduced a very high level of
precision in Lemma \ref{lem:d-cover}.  It is clear that the same
argument can replace the $1.03$ with $1.003$, etc.  These smaller
order errors are what accumulate into the decimal places in the final
$6.6$ which appears in Lemma \ref{lem:h-delete}, and which
later determines the value of $c$.  When $6.6$ is replaced by a
constant arbitrarily close to 6, the above argument produces a result
with $c$ approaching $1/13$.

\section{Concluding remarks}

Our proof of Theorem \ref{thm:3graph} is only valid when $n$ is
divisible by four, because we required a factor of 2 in each reduction
step (from hypergraphs to digraphs, and then to bipartite graphs).
Although we do not expect this condition to be necessary, removing
this restriction is open.  We also leave open the question of packing
Hamilton cycles of type $\ell$ in $k$-uniform hypergraphs when $k \geq
4$ and $\ell < k/2$.  Finally, another interesting direction is to
streamline the sets of pseudo-random properties which appear in the
statements of our packing results (both for digraphs and for
3-graphs).

\end{document}